\documentclass[twocolumn]{autart} 

\newtheorem{definition}{Definition}
\newtheorem{lemma}{Lemma}
\newtheorem{remark}{Remark}
\usepackage{bbm}
\usepackage{subfigure}
\usepackage{subfigmat}
\usepackage{amsfonts}
\usepackage{algorithm}
\usepackage{algpseudocode}
\usepackage{amsmath}
\usepackage{amssymb}
\usepackage{graphicx}      
\usepackage{xcolor}
\usepackage{comment}
\usepackage{tikz,pgfplots}
\usepackage{multirow}
\pgfplotsset{compat=1.15}
\usepackage{cite}
\usepackage{color,soul}

\usepackage{ulem} 

\newcommand{\XN}{\mathbb{X}_K}

\newcommand{\X}{\mathbb{X}}
\newcommand{\U}{\mathbb{U}}
\renewcommand{\Pr}[1]{\mathsf{Pr}\left\{#1\right\}}

\def\E#1{\mathbb{E}\left\{#1\right\}}

\newcommand{\R}{\mathbb{R}}

\def\R{\mathbb{R}}

\def\N{\mathbb{N}}

\def\tr#1{{\rm{tr}}\left(#1\right)}
\def\cov#1{{\rm{cov}}\left(#1\right)}

\everymath{\displaystyle}
\usepackage{xcolor}

\def\MF#1{\textcolor{black}{#1}}

\newtheorem{proposition}{Proposition}

\newtheorem {theorem}{Theorem}
\usepackage[utf8]{inputenc}

\def\Ec{{\mathcal{E}}}
\def\R{\mathbb{R}}

\def\N{\mathbb{N}}

\def\X{\mathbb{X}}

\def\z{\mathbf{z}}
\def\v{\mathbf{v}}
\def\g{\gamma}

\newcommand {\bsis} {\left\{ \begin{array} }
\newcommand {\esis} {\end{array}\right.}
 
\def\bmat#1{\left[\begin{array}{#1}}
\def\emat{\end{array}\right]}

\newcommand{\blista}{\renewcommand{\labelenumi}{(\roman{enumi})} 
	\begin{enumerate}}
	\newcommand{\elista}{\end{enumerate} \renewcommand{\labelenumi}{\arabic{enumi}.}}
\usepackage{color,soul}

\definecolor{mycolor1}{rgb}{0.00000,0.44700,0.74100}%
\definecolor{mycolor2}{rgb}{0.85000,0.32500,0.09800}%
\definecolor{mycolor3}{rgb}{0.92900,0.69400,0.12500}%

\usepackage{tikz}

\begin{document}

\begin{frontmatter}

\title{\MF{ Recursive feasibility for stochastic MPC \\ and the rationale behind fixing flat tires} }
%
\author[s1]{Mirko Fiacchini} \ead{mirko.fiacchini@gipsa-lab.fr}
\author[s2]{Martina Mammarella} \ead{martina.mammarella@cnr.it}
\author[s2]{Fabrizio Dabbene} \ead{fabrizio.dabbene@cnr.it}
\address[s1]{Univ. Grenoble Alpes, CNRS, Grenoble INP, GIPSA-lab, 38000 Grenoble, France.}
\address[s2]{Cnr-Istituto di Elettronica e di Ingegneria dell'Informazione e delle Telecomunicazioni, 10129 Torino, Italy.}

\begin{abstract}
In this paper, we address the problem of designing stochastic model predictive control (SMPC) schemes for linear systems affected by unbounded disturbances. The contribution of the paper is rooted in a measured-state initialization strategy.
First, due to the nonzero probability of violating chance-constraints in the case of unbounded noise, we introduce ellipsoidal-based probabilistic reachable sets and we include constraint relaxations to recover recursive feasibility conditioned to the measured state. Second, we prove that the solution of this novel SMPC scheme guarantees closed-loop chance constraints satisfaction under minimum relaxation. Last, we demonstrate that, in expectation, the need of relaxing the constraints vanishes over time, which leads the closed-loop trajectories steered towards the unconstrained LQR invariant region.
This novel SMPC scheme is proven to satisfy the recursive feasibility conditioned to the state realization, and its superiority with respect to open-loop initialization schemes  is shown through numerical examples.
\end{abstract}
\end{frontmatter}

\section{INTRODUCTION}
In real-world systems and control applications, safety and performance of the system may deteriorate due to several sources of uncertainty, typically paired with the complexity inherent to real world phenomena \cite{Prekopa1995}. In the framework of constrained dynamical systems, one may resort to \textit{robust} model predictive control (MPC) schemes, implicitly or explicitly addressing worst-case realizations of the disturbance \cite{mayne2005robust}. On the other hand, some conservativeness can be reduced whenever additional information about the uncertainty is available, e.g., in the form of probability distribution. In this case, one can rely on \textit{stochastic} MPC (SMPC) schemes, which have proved to be the state-of-the-art control approach for uncertain systems subject to constraints that are imposed in probability, i.e., formulated as chance constraints for which a certain amount of violation is permitted (see for instance \cite{farina2016stochastic,lorenzen2016constraint, munoz2020convergence} and references therein). 

Typically, SMPC solutions are classified as either \textit{randomized} methods, which rely on the generation of suitable disturbance realizations or scenarios \cite{Blackmore,CalafioreFagiano,lorenzen2016constraint}, or \textit{analytic approximation} methods, which exploit concentration inequalities -- as the classical Chebychev-Cantelli one -- to reformulate the probabilistic problem into a deterministic one (see \cite{farina2016stochastic, hewing2018stochastic} and references therein). Despite the class to which the different SMPC schemes belong, all share the common issue of the inherent difficulty of guaranteeing recursive feasibility of the underlying optimization problem. Indeed, this property, which represents a key feature of any MPC approach, {holding also in the robust context}, becomes very hard to {be ensured} in the probabilistic framework. 
For this reason, many SMPC approaches rely on the assumption of bounded disturbances \cite{kouvaritakis2016model,rawlings2017model} or bounded support distributions \cite{kouvaritakis2010explicit}. In this case indeed, recursive feasibility may be ensured by robust constraint tightening techniques, inspired by tube-based approaches. In particular, when the disturbance lies in a compact set, recursive feasibility can be recovered introducing a terminal cost and/or terminal constraints (see, e.g., \cite{goulart2006optimization}). 

On the other hand, there are many practical situations where this boundedness assumption is not realistic. In these cases, one has to account, by construction, for a \textit{nonzero probability} that the problem may become \textit{unfeasible}, since unbounded uncertainties almost surely lead to excursions of states from any bounded set \cite{paulson2020stochastic}. Consequently, in the SMPC framework, when facing unbounded uncertainties, guaranteeing closed-loop chance-constraints satisfaction and recursive feasibility may be difficult. The majority of approaches either rely on applying backup strategies in case of infeasibility, or 
on considering a ``fictitious state" (e.g., the predicted nominal at the previous instances) in spite of the measured one. In details, typical solutions rely on one of the two following strategies. The first class of approaches is based on the definition of a \textit{backup control scheme}, which is applied whenever the system states leave the region of attraction \cite{farina2013probabilistic,paulson2020stochastic}.
In the framework of backup control strategies, a natural choice to enable the state to be steered back to the region of attraction is to soften the state constraints, exploiting a strategy similar to the \textit{exact penalty function} method \cite{di1994exact}. Some example of backup control schemes have been proposed in \cite{yan2018stochastic}, where the chance constraints are defined as a discounted sum of violation probabilities along an infinite horizon, and in \cite{mammarella2020probabilistic}, where probabilistic validation techniques \cite{TeBaDa:97}, \cite{Alamo:15} have been used, combined with a penalty function method \cite{kerrigan2000soft,karg2019probabilistic}, to guarantee recursive feasibility without any assumption on independence or Gaussianity of the stochastic variables.

\MF{
The second class of approaches to handle recursive feasibility in the presence of unbounded disturbance relies instead on specific \textit{initialization strategies}, mainly alternating online between a closed-loop initialization, to be chosen as long as the problem is feasible, and a backup open-loop one, to be adopted whenever feasibility is lost for the current observed state \cite{farina2015approach,hewing2018stochastic}. However, this choice poses problems, since it purposely neglects the information carried by the current measurements whenever the measured states are not in the region of attraction of the controller. The first attempt to circumvent this issue is probably the one proposed in \cite{hewing2020recursively}, where an \textit{indirect feedback} over the measured state is introduced into the optimization problem through the cost function only, not in the constraints. Recursive feasibility is ensured by relying on probabilistic reachable sets (PRS) \cite{hewing2018stochastic,fiacchini2021probabilistic} to design tightened constraints for the nominal system, similar to tube-based MPC \cite{rawlings2017model}, provided the problem at the initial time is feasible. This research line was extended by introducing the idea of interpolation of the initial condition in \cite{kohler2022recursively} and \cite{schluter2022stochastic} and then generalized in  \cite{schluter2023stochastic}.
}

\MF{
This class of solutions, though, are based on the fact that, if the optimization problem is not feasible at the current state at some instant, then another optimization problem is solved at a fictitious state where feasibility holds. Such fictitious state can be either the state previously predicted or a convex combination between it and the real state. In any case, if constraints cannot be satisfied, the real state is disregarded for avoiding unfeasibility. The result is that the recursive feasibility and the probabilistic constraints satisfaction guarantee are holding {\it only with respect to the state measured at the initial time} (when the problem is assumed to be feasible), not for the real, measured one. In other words, the feasibility and constraints satisfaction of the open-loop solution computed at the initial time are the guarantees on which this approach rely.
But, of course, the feasibility for fictitious states is not an admissible solution, since neglecting the reality would not provide any reasonable insurance.
}

\MF{
To illustrate the concept with a common life situation, every initial probabilistic
insurance to not have a flat tire during the travel is useless once one has a puncture
along the way. The flat tire has to be fixed. For sure, one cannot act as if the puncture was not there and, of course, specific measures must be taken.
}

\MF{
The inherent contradiction of the methods proposed in the cited papers is highlighted also
by the fact that considering alternative initial conditions for the optimization problem, different from the real (measured) state, does not represent an admissible choice
at the initial instant. This leads to the absurd that a state that is unfeasible at the initial time might be feasible at later instants, by resorting to fictitious initialization strategies. In other words, the fact that the cited papers require initial
feasibility, implicitly recognizes that solving the problem for fictitious initial conditions is not an admissible choice.
}

\MF{
The underlying consideration on which our work is based is the realization of this
important fact: typically in SMPC, when facing unbounded uncertainties, guaranteeing closed-loop chance-constraints satisfaction and recursive feasibility may be impossible if one wants to rely on the measured state. Hence, one has to rely on some alternative strategies to restore those guarantees.}

\MF{The flat tire must be fixed.}

\MF{
In our paper, we propose an approach which differs from those available in the literature from a philosophical viewpoint. Indeed we claim that some kind of relaxation on the constraint sets is simply unavoidable, in general. This is specifically relevant if one does not want to neglect the information carried by the measured state. In fact, it is simple to observe that, in this framework, the future realization of the state could in general be arbitrarily far from the constraints set, in which case probabilistic constraints satisfaction might be impossible to hold. In this case, clearly, any probabilistic guarantee computed at the previous time steps becomes simply meaningless when a new realization of the state becomes available. Indeed, if the state at some point is far from the constraint sets, for instance, one cannot rely on the fact that some guarantees were valid in the previous instants. Since recursive feasibility consists of a guarantee that the deterministic optimization problem is solvable, there is no possibility of ensuring it in the future, but probabilistic guarantees on its feasibility can be given. 
}


The present work stems from the realization that, when dealing with unbounded stochastic uncertainty, and especially when a closed-loop initialization strategy is adopted, it becomes natural to introduce a concept of recursive feasibility using a \textit{probabilistic statement}. On the other hand, a probabilistic guarantee of feasibility is clearly not enough for our purposes, since we need to ensure that at \textit{every step} we are faced with a feasible problem, in order to be able to implement the approach. In particular, we propose a novel \textit{measured-state conditioned SMPC} (MS-SMPC) scheme where the recursive feasibility conditioned to the realization of the state is enabled by combining the adoption of probabilistic reachable sets with constraints relaxation. The underlying philosophy is rather simple. First, similar to \cite{farina2016stochastic, hewing2018stochastic}, we recast the stochastic optimization problem into a deterministic one, where, as in tube-based approaches, the probabilistic constraint sets are tightened through \textit{ellipsoidal-based PRS}, which size is directly related to the desired violation level. Then, we admit the constraints to be ``inflated'' by a factor $\gamma$ whenever the measured state is infeasible. This relaxation factor is introduced  into the cost function as penalty term, which allows us to employ a closed-loop initialization strategy, i.e.\  directly exploiting the current state measurement, while aiming to minimum relaxation. The idea of guaranteeing feasibility by relaxing the constraints has also been used in \cite{paulson2020stochastic} for treating the problem of hard input constraints in SMPC. However, unlike the approach in \cite{paulson2020stochastic}, in the setup we propose, the cost is linear in the constraints scaling parameters, and the overall problem is formulated as a convex (in particular, conic) optimization program. Additionally, we prove that, in expectation, the need of relaxation will vanish over time and the closed-loop trajectory will eventually converge to the LQR invariant region. 


The paper is structured as follows. Section~\ref{sec:SMPC_intro} introduces the novel SMPC scheme whereas the fundamental ingredients to compute the tightened constraint sets based on ellipsoidal PRS, thus providing the required probabilistic guarantees on chance constraint satisfaction, are presented in Section~\ref{sec:CC_approx}. The main technical results related to recursive feasibility conditioned to the measured state and the closed-loop constraint satisfaction under minimum relaxation conditions are detailed in Section~\ref{sec:main_res}. In Section~\ref{sec:add_res} we discuss the alternative approach of relaxing the probability of satisfaction (instead of enlarging the constraint sets), and different input initialization strategies. Last, the performance of the proposed {MS-SMPC} scheme are validated through numerical examples in Section~\ref{sec:num_sim}, where open- and closed- loop initialization approaches are compared. Main conclusions and future works are drawn in Section~\ref{sec:concl}.

{\small
\textbf{Notation.} The set $\mathbb{N}^{+}$ denotes the positive integers and $\mathbb{N} = \left\{0\right\} \cup\mathbb{N}^{+}$. 
Given $a,b\in\mathbb{N}$, $\mathbb{N}_a^{b}$ is the set of integers from $a$ to $b$. $A\ominus B=\left\{a\in A|\,a+b\in A, \forall b\in B\right\}$ denotes the Pontryagin set difference. We use $x_{\ell|k}$ for the state predicted $\ell$ steps ahead at time $k$, to differentiate it from the realization $x_{k+\ell}$. The sequence of length $N$ of vectors $v_{0|k}, \ldots, v_{N|k}$ is denoted in bold by $\textbf{v}_{k}$, its length begin left implicit. The expected value of a random variable $x$ is denoted $\E{x}$. With $W\succ 0$ ($W\succeq 0$) we denote a definite (semi-definite) positive matrix $W$. If $W\succeq 0$, then $W^{\frac{1}{2}}$ is the matrix satisfying $W^{\frac{1}{2}} W^{\frac{1}{2}}=W$. For $W\succeq 0$, we define $\|x\|_W\doteq \sqrt{x^\top W x}$. Given a matrix $W\succeq 0$ and a scalar $r \ge 0$, then $\Ec_{W}(r) =  \left\{ x\in\R^n\,|\, x^\top W^{-1} x \le r^2 \right\}$ represents the ellipsoid of ``{shape}'' $W$  and ``{radius}'' $r$. 
}

\section{PROBLEM FORMULATION}
\label{sec:SMPC_intro}
We consider a discrete-time, linear time-invariant system subject to additive disturbance
\begin{equation}\label{eq:sys}
x_{k+1}= A x_k + B u_k + w_k,
\end{equation}
with $x_k\in\mathbb{R}^n$, $u_k\in\mathbb{R}^m$, and $A,B$ of appropriate dimensions. The i.i.d. random noise $w_k \in \mathbb{R}^n$ is such that $\mathbb{E}\{w_k\} = 0$ and $\mathbb{E}\{w_k w_k^\top\}=\Gamma_w$ for all $k \in \N$. Note that, from independence, we have $\mathbb{E}\{w_i w_j^\top\} = 0$ for all $i \neq j$. 

\begin{remark}[On correlated noise]
\label{rem-correlation}
The latter assumption may be relaxed by assuming the existence of a \normalfont{correlation bound}\textit{, which requires only the existence of bounds on the mean and the covariance matrices, and a Schur stability condition on the closed-loop system, as proven in \cite{fiacchini2021probabilistic}.}
\end{remark}

Given the current state $x_k$, the states $x_{\ell|k}$ predicted $\ell$ steps ahead at $k$ obey the following dynamics
\begin{equation}\label{eq:xell}
x_{\ell+1|k}=Ax_{\ell|k}+Bu_{\ell|k}+w_{\ell+k},\quad x_{0|k}= x_k,
\end{equation}
where $x_{\ell|k}$ are random variables (due to the random noise assumption), and the inputs $u_{\ell|k}:\mathbb{R}^n\rightarrow\mathbb{R}^m$ are assumed to be measurable functions of $x_{\ell|k}$.

Due to the assumption of unbounded disturbances, no guarantee can be given on the recursive feasibility if one employs the measured state $x_k$, as it may inevitably lead to constraints violation. On the other hand, it would always be possible to relax the constraints on the state $\mathbb{X}$ and input $\mathbb{U}$, defined as  polytopic sets of the form
\begin{subequations}
\label{eq:XU}
\begin{align}
\mathbb{X} &\doteq \{ x \in \R^n\,|\, H_x x \leq h_x\}, \label{eq:X}\\
\mathbb{U} &\doteq \{ u \in \R^m\,|\, H_u u \leq h_u\}, \label{eq:U}
\end{align}
\end{subequations}
and containing the origin, by introducing scaling parameters $\gamma_x(x_k),\,\gamma_u(x_k)$ properly designed as functions of the realization $x_k$, to ensure the recursive feasibility and satisfaction of the closed-loop chance constraints, as shown later.

Hence, we assume that the system is subject to chance constraints on states and inputs with violation level\footnote{The choice of assuming the same $\varepsilon$ for both state and input is made to simplify the subsequent developments. Straightforwardly, the proposed approach can be extended to a more general framework accommodating diverse values.} $\varepsilon\in(0,1)$, where the probabilities are to be understood as conditioned on the measured state $x_k$, i.e.,
\begin{subequations}
\label{eq:constr_xk}
\begin{align}
&\Pr{x_{\ell|k}\in \gamma_x(x_{k})\mathbb{X} | \ x_k} \geq 1-\varepsilon,\label{eq:constr_xkx}\\
&\Pr{u_{\ell|k}\in \gamma_u(x_{k})\mathbb{U} | \ x_k} \geq 1 - \varepsilon.\label{eq:constr_xku}
\end{align} 
\end{subequations}
The dependence of $\gamma = ( \gamma_x, \gamma_u)$ on $x_k$ is left implicit, when clear from the context.  

\begin{remark}
 The proposed approach, for $\gamma_x=\gamma_u=1$, is analogous to the one employed in \cite{hewing2018stochastic,hewing2020recursively} with the following main differences: i) initial condition set as $z_{0|k} = z_{1|k-1}$, ii) chance constraints defined as
\begin{align*}
\Pr{x_k \in \mathbb{X} | x_0} \geq 1-\varepsilon_x, \quad \Pr{u_k \in \mathbb{U} | x_0} \geq 1 - \varepsilon_u, 
\end{align*}
which satisfaction for all $k \in \N^+$ is conditioned to the feasibility of the optimization problem for $x_0$ at $k=0$.
\end{remark}

Then, as typically done in the SMPC framework, we consider as terminal set a probabilistic invariant set $\mathbb{X}_N$ contained in 
\begin{equation}
\label{eq:XN}
    \XN \doteq\Big\{x\in \mathbb{X}\,|\, 
    H_u K x\leq h_u
    \Big\},
\end{equation}
i.e., the set of states satisfying both state and input constraints under the stabilizing feedback law $u_k=Kx_k$. Specifically, we consider the LQR state-feedback control law $u_{k}=Kx_{k}$,
with $K\in\mathbb{R}^{m\times n}$, i.e., given $Q\succeq 0$, $R\succ 0$, a Lyapunov matrix $P\succ 0$ exists satisfying
\begin{equation}
\label{eq:LQR}
Q+K^{\top}RK+A_K^\top PA_K-P =  0
\end{equation}
thus guaranteeing that $A_{K} \doteq A+BK$ is Schur.

Moreover, as typical in stabilizing MPC,  we introduce the terminal cost $V_f(x_{N|k})=\|x_{N|k}\|_P^2$ with $P\succ0$ the solution of \eqref{eq:LQR}, so that the finite horizon cost $J_N(\textbf{x}_k,\textbf{u}_k)$ to be minimized at time $k$ is defined as
\begin{equation}
    \label{eq:init_cost}
    J_N(\mathbf{x}_k,\mathbf{u}_k) \doteq\mathbb{E}\Bigg\{\sum_{\ell=0}^{N-1}\Big( \|x_{\ell|k}\|_Q^2+\|u_{\ell|k}\|^2_R\Big)+V_f(x_{N|k})\Bigg\}.
\end{equation}

Hence, the resulting stochastic (finite-horizon) optimal control problem can be stated as
\begin{equation}
\label{eq:SMPC_init_S}
\begin{aligned}
& \hspace{-0.7cm} \min_{\substack{\textbf{x}_k, \textbf{u}_k\\ \gamma_x, \gamma_u,\gamma_N}}  \mathbb{E}\{J_N(\mathbf{x}_k,\mathbf{u}_k)\}\\
 \text{s.t.} \ & x_{\ell+1|k} = A x_{\ell|k} + B u_{\ell|k} + w_{k+\ell},\\
& x_{0|k}={x_k},\\  
& \Pr{x_{\ell|k}\in \gamma_x\mathbb{X}\,\,| \ x_k } \geq 1-\varepsilon,\quad \ell\in\mathbb{N}_1^{N-1}\\
& \Pr{u_{\ell|k}\in \gamma_u\mathbb{U}\,\,| \ x_k } \geq 1 - \varepsilon,\quad \ell\in\mathbb{N}_1^{N-1}\\
& \Pr{x_{N|k}\in\gamma_N\mathbb{X}_N\,\,| \ x_k }\geq 1-\varepsilon.
\end{aligned}
\end{equation}
%
%
To compute explicitly the expected value of $J_N(\mathbf{x}_k,\mathbf{u}_k)$ in \eqref{eq:init_cost}, we recast the chance-constrained optimization problem \eqref{eq:SMPC_init_S} into a deterministic one resorting to a tube-based approach inherited from robust MPC schemes (see for instance \cite{mammarella2018offline}). In this case, the predicted state $x_{\ell|k}$ is split into a deterministic, nominal part $z_{\ell|k}$, and a stochastic, error part $e_{\ell|k}$ such that
\begin{equation}
    x_{\ell|k}=z_{\ell|k}+e_{\ell|k}.
    \label{eq:nom_err}
\end{equation}
Then, we consider a prestabilizing error feedback, which leads to the following predicted input
\begin{equation}\label{eq:input}
    u_{\ell|k}=Ke_{\ell|k}+v_{\ell|k},
\end{equation}
where $K$ the solution of \eqref{eq:LQR}, and $v_{\ell|k}$ are the free SMPC optimization variables. Hence, the nominal system and error dynamics are given, {respectively}, by
\begin{subequations}
\label{eq:sys_ze}
\begin{align}
z_{\ell+1|k} &= A z_{\ell|k} + B v_{\ell|k}, \quad z_{0|k}=z_{k}, \label{eq:sys_z}\\
e_{\ell+1|k} &= A_Ke_{\ell|k} + w_{k+\ell},\quad e_{0|k}=0, \label{eq:sys_e}
\end{align}
\end{subequations}
where $e_{\ell|k}$ are zero-mean, random variables.

Correspondingly, it is possible to rely on \eqref{eq:sys_ze} to redefine $J_N(\mathbf{x}_k,\mathbf{u}_k)$ thus obtaining
\begin{align*}
        J_N(\mathbf{x}_k,\mathbf{u}_k) 
        &=\sum_{\ell=0}^{N-1}\Big( \|z_{\ell|k}\|_Q^2+\|v_{\ell|k}\|^2_R\Big)+\|z_{N|k}\|_P^2+c,
\end{align*}
$$c=\mathbb{E}\Bigg\{\sum_{\ell=0}^{N-1}\Big(\|e_{\ell|k}\|_{Q+K^\top RK}^2\Big)+\|e_{N|k}\|_P^2\Bigg\}$$
where $c$ is a constant term and can be neglected in the optimization. In this way, we obtain a quadratic, finite-horizon cost
\begin{equation}\label{eq:Jp}
J_p(\mathbf{z}_k,\mathbf{v}_k) =\sum_{\ell=0}^{N-1}\Big(\|z_{\ell|k}\|_Q^2+\|v_{\ell|k}\|^2_R\Big)+\|z_{N|k}\|_P^2,
\end{equation}
depending on the deterministic variables $z_{\ell|k}$ and $v_{\ell|k}$.

Then, we define a term, to be added to the cost function $J_p(\mathbf{z}_k,\mathbf{v}_k)$, that accounts for the constraints relaxation required to ensure recursive feasibility, i.e.,
\begin{equation}\label{eq:Jr}
J_r(\gamma)\doteq \max\{\gamma_x-1,\gamma_u-1\},
\end{equation}%
to minimize the constraint relaxation for a given $x_k$. Hence, the deterministic cost function to be minimized is given by
\begin{align}
J_N(\mathbf{z}_k,\mathbf{v}_k,\gamma)& = J_p(\mathbf{z}_k,\mathbf{v}_k)+\eta J_r(\gamma),\label{eq:cost_z_tot}
\end{align}
which aims to balance performance and constraints relaxation to ensure recursive feasibility when the deterministic optimization problem is initialized with the measured state $x_k = z_{0|k}$.

Finally, we can recast the chance-constrained probabilistic optimization problem \eqref{eq:SMPC_init_S} into a novel deterministic, measured-state dependent, and convex optimization problem of the form
\begin{subequations}
\label{eq:SMPC_mixedn}
\begin{align}
&\hspace{-0.3cm}  \min_{\textbf{z}_k, \textbf{v}_k, \gamma_x, 
\gamma_u} 
J_p(\mathbf{z}_k,\mathbf{v}_k)+ \eta J_r(\gamma)\label{eq:cost-balancedn}\\
&\hspace{-0.3cm} \text{s.t.} \ z_{\ell+1|k} = A z_{\ell|k} + B v_{\ell|k},\\ 
&\hspace{-0.3cm} \phantom{s.t.} \ z_{0|k} = x_k,\\ 
&\hspace{-0.3cm} \phantom{s.t.} \ z_{\ell|k} \in \gamma_x\Ec_{W_x}(r_x)\ominus \mathcal{R}_\ell^x(\varepsilon), \; \ell\in\mathbb{N}_1^{N-1}, \label{eq:prabarBn}\\
&\hspace{-0.3cm} \phantom{s.t.} \ v_{\ell|k} \in \gamma_u\Ec_{W_u}(r_u)\ominus \mathcal{R}_\ell^u(\varepsilon), \; \ell\in\mathbb{N}_1^{N-1}, \label{eq:prbbarBn}\\
&\hspace{-0.3cm} \phantom{s.t.} \ z_{N|k} \in \gamma_x\Ec_{W_x}(r_x)\ominus \mathcal{R}_N(\varepsilon), \label{eq:prcbarBn}\\
&\hspace{-0.3cm} \phantom{s.t.}  \ z_{N|k} \in \gamma_u\Ec_{W_x}(r_u)\ominus \mathcal{R}_N(\varepsilon), \label{eq:prcbar2Bn}\\
&\hspace{-0.3cm} \phantom{s.t.} \ \gamma_x\ge 1, \ \gamma_u\ge 1 \label{eq:constrgamman}
\end{align}
\end{subequations}
where $\Ec_{W_x}(r_x)$ and $\Ec_{W_u}(r_u)$ are ellipsoidal approximations of $\X$ and $\U$, and $\mathcal{R}_\ell^x(\varepsilon)$, $\mathcal{R}_\ell^u(\varepsilon)$, and $\mathcal{R}_N(\varepsilon)$ are probabilistic reachable sets (PRS) properly designed to enforce the required probabilistic properties to the solution. Last, in Section \ref{sec:inputbounds} we will outline three possible choices for the initial bound on $v_{0|k}$, each one providing a specific feature to the optimization problem: (i) \textit{no bounds}, (ii) \textit{hard bounds}, and (iii) \textit{soft bounds}.

\begin{remark}
Note that, although the cost is not differentiable, the problem \eqref{eq:SMPC_mixedn} is still convex. Moreover, it is possible to consider smooth approximations, e.g., by replacing the $\max$ in \eqref{eq:SMPC_mixedn} with $\gamma_x^{\beta} + \gamma_u^{\beta}$, with $\beta \geq 1$.
\end{remark}

In the following section, we will present the approach to compute the tightened sets \eqref{eq:prabarBn}-\eqref{eq:prcbar2Bn} based on ellipsoidal PRS to attain the required probabilistic guarantees for Problem~\eqref{eq:SMPC_mixedn}, conditioned  to the realization $x_k$. Specifically,  under suitable conditions and properly designing the PRS, we will prove that:
\begin{enumerate}
\item[1.] the solution to Problem~\eqref{eq:SMPC_mixedn} guarantees \textit{closed-loop constraints satisfaction under minimum relaxation} and \textit{recursive feasibility}, conditioned to $x_k$.
\item[2.] \textit{the need of relaxation vanishes} in expectation while the closed-loop trajectories are steered towards a specific invariant region.
\end{enumerate}

\section{Chance constraints ellipsoidal tightening}
\label{sec:CC_approx}
The concept of probabilistic reachable sets in the framework of SMPC was used in several setup as \cite{hewing2018stochastic}. As described in \cite{hewing2020recursively}, PRS are employed for constraint tightening, playing a role similar to robust invariant sets in tube-based MPC. In the following, we extend this concept to the framework of ellipsoidal-based PRS for the nominal system \eqref{eq:sys_z}, starting from computing PRS for the error system \eqref{eq:sys_e}, such that probabilistic guarantees are achieved in terms of recursive feasibility and chance-constraint satisfaction. Moreover, we demonstrate how the proposed design of PRS-based tightened constraint sets allows to reduce some conservativeness of other approaches relying on invariant sets for the nominal dynamics (see Remark \ref{rem:conserv}).   

\subsection{Preliminary definitions}
First, we report the definitions of probabilistic reachable and invariant sets generally employed in the SMPC framework, e.g., see \cite{kofman2012probabilistic,hewing2018correspondence,fiacchini2021probabilistic}.

\begin{definition}[Probabilistic reachable set] The sequence $\mathcal{R}_\ell$, $\ell\in\mathbb{N}$, is a sequence of PRS for the system $\xi_{\ell+1} = A\xi_\ell+w_\ell$, with violation level $\varepsilon\in[0,1)$, if 
$\xi_0\in\mathcal{R}_0$ implies $\Pr{\xi_\ell\in\mathcal{R}_\ell} \geq 1-\varepsilon $ for all $\ell\in\mathbb{N}^{+}$.
\end{definition}
\begin{definition}[Probabilistic invariant set]
The set $\Omega \subseteq \R^n$ is a PIS for the system 
$\xi_{\ell+1} = A\xi_\ell+w_\ell$, with violation level $\varepsilon \in [0,1)$,  if $ \xi_0 \in \Omega$ implies $\Pr{\xi_\ell \in \Omega} \geq 1-\varepsilon $ for all $\ell\in\mathbb{N}^{+}$.
\end{definition}

Next, we present the application of this definition to the design of ellipsoidal PRS to be enforce in Problem~\eqref{eq:SMPC_mixedn} to achieve the required probabilistic guarantees.

\subsection{Design of ellipsoidal PRS}
Let us consider the error system \eqref{eq:sys_e}. Following the assumptions on the disturbance $w_{k+\ell}$, we have that the mean is $\E{e_{\ell|k}} = 0$ for all $\ell$, whereas the covariance matrix $E_{\ell|k}\doteq \mathbb{E}\{e_{\ell|k} e_{\ell|k}^\top\}$ satisfies the following recursion
\begin{equation}
E_{\ell+1|k} = A_K E_{\ell|k} A_K^\top + \Gamma_w, \label{eq:Ek}
\end{equation}
for all $k\in\mathbb{N}$, with $E_{0|k} = 0$. Then, as discussed in \cite{kofman2012probabilistic,hewing2018correspondence,fiacchini2021probabilistic},  {the Chebychev inequality with the covariance matrix recursion provided in \eqref{eq:Ek}  guarantees that}
\begin{equation}
\label{eq:prob_phi}
{\Pr{e_{\ell|k}\in\Ec_{E_{\ell|k}}\!(\rho)}} \geq 1- \varepsilon,
\end{equation}
where $\rho$ is a strictly decreasing function of~$\varepsilon$, i.e., 
\begin{equation}
    \label{eq:eps_r}
\rho\,:\,
\begin{cases}
   \frac{n}{\rho^2} = \varepsilon &\text{ for generic distributions},\\
  1 - \chi_n^2(\rho^2)  = \varepsilon &\text{ for Normal distributions}.
\end{cases}
\end{equation}
Now, we exploit the previous result to construct a sequence of PRS for the error system \eqref{eq:sys_e} subject to \textit{uncorrelated} noise, relying on an approach similar to the one in \cite[Proposition 4]{fiacchini2021probabilistic}. Specifically, every matrix $W_x\succ 0$ and parameter $\lambda \in (0, 1)$ satisfying
\begin{subequations}
\label{eq:reachLMI}
    \begin{align}
 A_K W_x A_K^\top &\preceq \lambda^2 W_x  , \label{eq:reachLMI_1}\\
 \Gamma_w  &\preceq (1 - \lambda)^2 W_x, \label{eq:reachLMI_2}
\end{align}
\end{subequations}
are such that, for $\rho>0$, the sequence of ellipsoids of shape $W_x$ and increasing radius $\rho( 1-\lambda^\ell)$, i.e.,
\begin{equation}\label{eq:PRS_ell}
\mathcal{R}_\ell^x(\varepsilon) = \Ec_{W_x}\!(\rho(1-\lambda^\ell)),\quad {\ell\in\mathbb{N}},
\end{equation} 
is a sequence of PRS for \eqref{eq:sys_e} with violation probability $\varepsilon$. Formally, given $\rho$ as a function of $\varepsilon$ \eqref{eq:eps_r}, we have 
\begin{equation}\label{eq:Pre}
    \Pr {e_{\ell|k}\in \Ec_{W_x}\!(\rho(1-\lambda^\ell))}\geq 1-\varepsilon \quad \forall \ell\in\mathbb{N}.
\end{equation}
Once selected the shape $W_x$ of the PRS for the error system \eqref{eq:sys_e}, we design the maximum-volume ellipsoid inscribed in $\mathbb{X}$ of the same shape and radius $r_x$, i.e., 
\begin{equation}\label{eq:Wx}
\Ec_{W_x}\!(r_x) = \{x \in \R^{n}: \ x^\top W_x^{-1} x \leq r_x^2\} \subseteq \X.
\end{equation}
\begin{remark}
Note that the possible choice of parameters $K$, $W_x$, $\lambda$, and $r_x$ is not unique. One, in fact, may want to obtain a set $\Ec_{W_x}\!(r_x)$ providing a better approximation of $\X$. Alternatively, one may select a small $\lambda$ to achieve faster convergence for the nominal state, which  might lead though to an aggressive local control $K$ and then to very tight input constraints on $v_{\ell|k}$.  
\end{remark}
In the following proposition, we prove that if we design the tightened state constrain set as
$$\mathbb{Z}_\ell\doteq \Ec_{W_x}(r_x) \ominus \mathcal{R}_{\ell}^x(\varepsilon),$$
then the state chance constraints conditioned to $x_k$, i.e.,
$\Pr{x_{\ell|k}\in \mathbb{X} | \ x_k} \geq 1 - \varepsilon$,
are satisfied for all $\ell\in\N$.

\begin{proposition}[Tightened state constraints]
\label{prop:xellk}
Consider the system \eqref{eq:sys_ze}, $W_x\succ 0$ and $\lambda \in (0,1)$ satisfying \eqref{eq:reachLMI}, and $r_x$ as in \eqref{eq:Wx}.
For each $\ell \in \N$, let $\rho>0$ be such that $\rho(1-\lambda^\ell)\le r_x$. If $z_{\ell|k}$  fulfills
\begin{align}
\label{eq:prop1}
z_{\ell|k}\in\mathbb{Z}_\ell&=\Ec_{W_x}(r_x)\ominus\Ec_{W_x}\big(\rho(1-\lambda^\ell)\big)\nonumber\\
& =\Ec_{W_x}\!\left(r_x-\rho(1-\lambda^\ell) \right)
\end{align}
then the state chance-constraint conditioned to $x_k$, i.e.,
\begin{equation}\label{eq:Prop1a}
    \Pr{x_{\ell|k}\in \mathbb{X} \ | \ x_k} \geq 1- \varepsilon, 
\end{equation} 
is satisfied. 
\end{proposition}
The proof is reported in Appendix~\ref{app:prop1}.

Now, we proceed with an analogous approach to design the tightened input constraint set $\mathbb{V}_\ell$. First, we compute the PRS $\mathcal{R}_\ell^u$ in \eqref{eq:prbbarBn}, which is the set of $e_{\ell|k}$ such that $Ke_{\ell|k}\in\mathcal{R}_\ell^u$. This PRS is defined as
\begin{equation}\label{eq:PRS_ell_u}
\mathcal{R}_\ell^u(\varepsilon) = \Ec_{W_u}(\rho (1-\lambda^\ell)),  \quad \forall \ell \in \N,
\end{equation} 
where the shape matrix $W_u$ satisfies
\begin{equation}\label{eq:WuWx}
K^\top W_u^{-1} K \preceq W_x^{-1}.
\end{equation}
Then, given the matrix $W_u$, the radius $r_u$ is computed as  
\begin{equation}\label{eq:Wu}
\Ec_{W_u}(r_u) = \{u \in \R^{m}: \ u^\top W_u^{-1} u \leq r_u^2\} \subseteq \U,
\end{equation}
such that $\Ec_{W_u}(r_u)$ is the maximal-volume ellipsoid inscribed in $\U$. In the following proposition, which proof is reported in Appendix~\ref{app:prop2}, we demonstrate that if one defines $\mathbb{V}_\ell=\Ec_{W_u}(r_u)\ominus\mathcal{R}_\ell^u(\varepsilon)$, the input chance constraint conditioned to $x_k$, i.e.,
$\Pr{u_{\ell|k}\in \mathbb{U} | \ x_k} \geq 1 - \varepsilon$, is satisfied, provided that $v_{\ell|k} \in \mathbb{V}_\ell$.

\begin{proposition}[Tightened input constraints]
\label{prop:uellk}
Given the system \eqref{eq:sys_ze}, $W_x\succ 0$ and $\lambda \in (0,1)$ satisfying \eqref{eq:reachLMI}, and $W_u\succ 0$ and $r_u>0$ such that \eqref{eq:WuWx} and \eqref{eq:Wu} hold. For each $\ell \in \N$, let $\rho>0$ be such that $\rho (1-\lambda^\ell) \le r_u$. If $v_{\ell|k}$
fulfills
\begin{align}
\label{eq:prop2}
v_{\ell|k}\in\mathbb{V}_\ell&=\Ec_{W_u}(r_u)\ominus\Ec_{W_u}\big(\rho(1-\lambda^\ell)\big)\nonumber\\
& =\Ec_{W_u}\!\big(r_u-\rho (1-\lambda^\ell)\big)
\end{align}
then the input chance-constraint conditioned to $x_k$, i.e.,
\begin{equation}
    \Pr{u_{\ell|k}\in \mathbb{U} | \, x_k} \geq 1 - \varepsilon,
\end{equation}
is satisfied.  
\end{proposition}


Next, we provides the main features to design the ellipsoidal approximation of the terminal constraint set $\XN$.
\begin{lemma}\label{lem:rW}
Let matrices $W_x, W_u$ and radii $r_x, r_u$ satisfy \eqref{eq:Wx}, \eqref{eq:WuWx}, and \eqref{eq:Wu}. The ellipsoidal approximation of shape $W_x$ of the set $\X_K$ is given by
\begin{equation}\label{eq:RN}
\Ec_{W_x}(r_{xu}) \subseteq \XN,
\quad r_{xu} \doteq \min\{r_x, \, r_u\},  
\end{equation}
and it is such that, if $x \in \Ec_{W_x}(r_{xu})$, then $x \in \Ec_{W_x}(r_{xu}) \subseteq \X$ and $u = Kx \in \Ec_{W_u}(r_{xu}) \subseteq \U$.  
\end{lemma}
The proof to Lemma~\ref{lem:rW} is reported in Appendix~\ref{app:lemma1}.

Then, to design the tightened terminal constraint set $\mathbb{Z}_N$, it is necessary to search for a condition on the nominal state at time $\ell = N$ that implies the chance constraints satisfaction at time~$N$ and also in the whole future. 
To do this, we rely on Lemma~\ref{lem:rW} and on the concepts of PRS and PIS, as shown in the next theorem, which proof can be found in Appendix~\ref{app:th_term}.

\begin{theorem}\label{th:term}
Given the system \eqref{eq:sys_ze}, $W_x \succ 0$ and $\lambda \in (0,1)$ satisfying \eqref{eq:reachLMI}, and $r_{xu}= \min\{r_x, r_u\}$ such that \eqref{eq:Wx}, \eqref{eq:Wu}  and \eqref{eq:WuWx} hold, and $\rho \leq r_{xu}$. If $z_{N|k}$ satisfies 
\begin{align}
\label{eq:constr_zNk}
z_{N|k}\in\mathbb{Z}_N&=\Ec_{W_x}\!(r_{xu})\ominus \mathcal{R}_N(\varepsilon)\nonumber\\
& =\Ec_{W_x}\!\left(r_{xu} - \rho(1-\lambda^N) \right)
\end{align}
then the terminal chance constraint 
\begin{align}
&\Pr{x_{\ell|k}\in \XN | x_k} \geq 1-\varepsilon,\label{eq:constr_xkXN}
\end{align} 
with $v_{\ell|k} = K z_{\ell|k}$ is satisfied for all $\ell = N + j$ with $j \in \N$. 
\end{theorem}

\begin{remark}
In Theorem \ref{th:term}, the terminal region for the nominal system is defined with respect to the radius $r_{xu}=\min\{r_x,r_u\}$. This is required to guarantee that both $z_{N|k}\in\Ec_{W_x}(r_x)\ominus\mathcal{R}_N(\varepsilon)$ and $z_{N|k}\in\Ec_{W_x}(r_u)\ominus\mathcal{R}_N(\varepsilon)$ are simultaneously satisfied at step N.  
\end{remark}
Additionally, it is important to highlight that if the set
$$\lim_{j \rightarrow +\infty} \Ec_{W_x}\!(\rho(1-\lambda^j)) = \Ec_{W_x}\!(\rho),$$
which is an outer approximation of the minimal probabilistic invariant set with violation probability $\varepsilon$, is not contained in $\mathbb{X}$, then the chance constraint on the state may be violated along the trajectory. This gives a geometric meaning to constraint $\rho \leq r_{xu}$.

\begin{remark}
\label{rem:conserv}
It is important to remark that the classical approach to guarantee recursive satisfaction of the terminal constraint would lead to a solution more conservative than the one proposed in this paper. Indeed, a common choice consists in imposing that the final state shall belong to an invariant set for the nominal dynamics. In the ellipsoidal setting, this could be done by imposing 
\begin{equation}\label{eq:term_inv}
z_{N|k} \in \Ec_{W_x}\!(r_{xu}) \ominus \Ec_{W_x}\!(\rho) = \Ec_{W_x}\!(r_{xu} - \rho).
\end{equation}
It is evident, by comparing \eqref{eq:constr_zNk} with \eqref{eq:term_inv}, that the proposed terminal region is less conservative than the one obtained with classical approaches relying on the concept of terminal invariance.
\end{remark}
\subsection{MS-SMPC algorithm}
Before presenting the algorithm associated to the proposed MS-SMPC, it is worth highlight that any relaxation on the ellipsoids can be seen as a relaxation on the state and input polytopes.
\begin{remark}[Scaling of  ellipsoid approximations]\label{rem:scaling}
Any scaling of the sets $\X$, $\U$, and $\XN$ would lead to the same scaling for their inner approximating ellipsoids $\Ec_{W_x}(r_x)$, $\Ec_{W_u}(r_u)$, and $\Ec_{W_x}(r_{xu})$. Indeed, for any generic scaling parameter $\gamma\geq 0$, we have
\begin{align*}
&\Ec_{W_x}(\gamma r_x) = \gamma \Ec_{W_x}(r_x) \subseteq \gamma \X,\\
&\Ec_{W_u}(\gamma r_u) = \gamma \Ec_{W_u}(r_u) \subseteq \gamma \U,\\
&\Ec_{W_x}(\gamma r_{xu}) = \gamma \Ec_{W_x}(r_{xu}) \subseteq \gamma \XN.
\end{align*}
\end{remark}
Hence, according to the specific definition of the tightened sets introduced above and Remark~\ref{rem:scaling}, the MPC optimization problem \eqref{eq:SMPC_mixedn} takes the form of the following conic programming problem   
\begin{subequations}\label{eq:SMPC_mixed}
\begin{align}
&  \min_{\z_k, \v_k, \g_x, 
\g_u} \! \! \sum_{\ell = 0}^{N-1} \left(\|z_{\ell|k}\|_Q^2+ \|v_{\ell|k}\|_R^2\right) + \|z_{N|k}\|_P^2) \nonumber\\ 
& \hspace{2.3cm} +\,\,\eta \, \max\{\gamma_x-1, \, \gamma_u-1\} \label{eq:cost-balanced}\\
&\hspace{-0.cm} \text{s.t.} \ z_{\ell+1|k} = A z_{\ell|k} + B v_{\ell|k},\\ 
&\hspace{-0.cm} \phantom{s.t.} \ z_{0|k} = x_k,\\ 
&\hspace{-0.cm} \phantom{s.t.} \ z_{\ell|k}^\top W_x^{-1} z_{\ell|k} \leq (\gamma_x r_x - \rho(1-\lambda^\ell))^2, \; \ell\in\mathbb{N}_1^{N-1}, \label{eq:prabar}\\
&\hspace{-0.cm} \phantom{s.t.} \ v_{\ell|k}^\top W_u^{-1} v_{\ell|k} \leq (\gamma_u r_u - \rho(1-\lambda^\ell))^2, \; \ell\in\mathbb{N}_1^{N-1}, \label{eq:prbbar}\\
&\hspace{-0.cm} \phantom{s.t.} \ z_{N|k}^\top W_x^{-1} z_{N|k} \leq (\gamma_x r_x - \rho(1-\lambda^N))^2, \label{eq:prcbar}\\
&\hspace{-0.cm} \phantom{s.t.} \ z_{N|k}^\top W_x^{-1} z_{N|k} \leq (\gamma_u r_u - \rho(1-\lambda^N))^2, \label{eq:prcbar2}\\
&\hspace{-0.cm} \phantom{s.t.} \ 1 \leq \gamma_x, \quad 1 \leq \gamma_u, \label{eq:constrgamma}
\end{align}
\end{subequations}
to be solved at every sample instant to compute the control input to be applied $u_k^*(x_k) = v_{0|k}^*$.

Finally, we present the algorithm associated with the proposed MS-SMPC, defined by Problem~\eqref{eq:SMPC_mixed}.

\noindent\textsc{Offline Step.}
\begin{enumerate}{\it
\item[1.] Select $W_x\succ0$ and $\lambda\in(0,1)$ satisfying \eqref{eq:reachLMI}.
\item[2.] Compute $r_x$ as in \eqref{eq:Wx} to obtain $\Ec_{W_x}(r_x)\subseteq\X$.
\item[3.] Compute $W_u\succ 0$ and $r_u>0$ such that \eqref{eq:WuWx} and \eqref{eq:Wu} hold, obtaining $\Ec_{W_u}(r_u)\subseteq\U$.
\item[4.] Select $\rho$ according to $\varepsilon\in(0,1)$ as in \eqref{eq:eps_r}.
\item[5.] Define the PRS $\mathcal{R}_\ell^x(\varepsilon)$, $\mathcal{R}_\ell^u(\varepsilon)$, and $\mathcal{R}_N(\varepsilon)$ according to \eqref{eq:PRS_ell}, \eqref{eq:PRS_ell_u}, and \eqref{eq:RN}, respectively.  
} 
\end{enumerate}
\textsc{Online Implementation.}  {\it At each time $k$:}
\begin{enumerate}{\it
\item[1.] Measure the current state $x_{k}$.
\item[2.] Solve the convex optimization problem \eqref{eq:SMPC_mixed}.}
\it{\item[3.] Apply the control input} $u_{k}=v^\star_{0|k},$ \it{where $v_{0|k}^\star$ is the first element of the optimal control sequence} $\mathbf{v}_{k}^\star$.
\end{enumerate}
\section{Properties of MS-SMPC}
\label{sec:main_res}
In this section, we prove how the solution to Problem~\eqref{eq:SMPC_mixed} guarantees recursive feasibility conditioned to the measured state $x_k$ and closed-loop chance constraints satisfaction under minimum relaxation, thanks to the combination of ellipsoidal PRS and constraints scaling. 

First, a preliminary result is introduced which demonstrates that, for $\gamma_x^\star = \gamma_u^\star = 1$, the optimal solution of Problem~\eqref{eq:SMPC_mixed} is the LQR control.
\begin{proposition}\label{prop:inEcWx}
If $x_k \in \Ec_{W_x}\!(r_{xu})$, then the optimal solution to Problem~\eqref{eq:SMPC_mixed} is given by 
$z^\star_{\ell+1|k} = A_K z^\star_{\ell|k}$ for all $\ell \in \N_{0}^{N-1}$, with $v^\star_{\ell|k} = K z^\star_{\ell|k}$, $\gamma_x^\star = 1$, and $\gamma_u^\star = 1$. Then, $u^\star_{0|k} = K x_k$ and $J_N^\star(x_k) = J_{LQR}(x_k) = x_k^\top P x_k$.
\end{proposition}
The proof is reported in Appendix~\ref{app:Prop3}.

Given the optimal solution $\z^\star(x_k),\,\v^\star(x_k)$, $\gamma_x^\star(x_k)$, and $\gamma_u^\star(x_k)$ to Problem~\eqref{eq:SMPC_mixed}, we define the candidate solution at time $k+1$ as ${\bf z}^+(x_k,w_k) = \{z_{\ell|k+1}\}_{\ell=0}^N$ and ${\bf v}^+(x_k,w_k) = \{v_{\ell|k+1}\}_{\ell=0}^{N-1}$ with
\begin{subequations}\label{eq:zvplus}
\begin{align}
& z_{\ell|k+1} = z_{\ell+1|k}^\star + A_K^{\ell} w_k, \qquad \forall \ell \in \N_0^{N-1},\\
& z_{N|k+1} = A_K z_{N|k}^\star + A_K^{N} w_k,\\
& v_{\ell|k+1} = v_{\ell+1|k}^\star + K A_K^{\ell} w_k, \quad \forall \ell \in \N_0^{N-2},\label{eq:th1v}\\
& v_{N-1|k+1} = K z_{N|k}^\star + K A_K^{N-1} w_k. \label{eq:th1vN}
\end{align}
\end{subequations}
Hence, the candidate solutions at time $k+1$ depends on the realization of $x_k$ and $w_k$. To simplify the notation, in what follows the explicit dependencies on $x_k$ and $w_k$ will be omitted when clear from the context.

The following theorem provides a candidate feasible solution for the relaxation parameters at time $k+1$, given the optimal solution at time $k$, i.e., $\z^\star(x_k), \v^\star(x_k), \g^\star_x(x_k), \g^\star_u(x_k)$, and a candidate nominal trajectory $\z^+(x_k,w_k), \v^+(x_k,w_k)$ at time $k+1$ given by \eqref{eq:zvplus}. 
\begin{theorem}\label{th:rplus}
Consider the optimal solution to Problem~\eqref{eq:SMPC_mixed} at time $k$, the realization of the disturbance $w_k$, and the candidate solution at time $k+1$ given by ${\bf z}^+(x_k,w_k)$ and ${\bf v}^+(x_k,w_k)$ defined in \eqref{eq:zvplus}. Define the candidate scaling parameters at time $k+1$ as
\begin{subequations}
\label{eq:rxplus}
\begin{align*}
& \hspace{-0.cm} \gamma_x^+(x_k,w_k) = \arg\min_{\gamma_x} \ \gamma_x\\ 
&\hspace{-0.cm} \text{s.t.} \ z_{\ell|k+1}^\top W_x^{-1} z_{\ell|k+1} \leq (\gamma_x r_x - \rho(1-\lambda^\ell))^2, \; \ell\in\mathbb{N}_1^{N-1}, \\
&\hspace{-0.cm} \phantom{s.t.} \ z_{N|k+1}^\top W_x^{-1} z_{N|k+1} \leq (\gamma_x r_x - \rho(1-\lambda^N))^2, \\
&\hspace{-0.cm} \phantom{s.t.} \ 1 \leq \gamma_x,\\
& \textrm{ and } \\
& \hspace{-0.cm} \gamma_u^+(x_k,w_k) = \arg\min_{\gamma_u} \ \gamma_u\\ 
&\hspace{-0.cm} \text{s.t.} \ v_{\ell|k+1}^\top W_x^{-1} v_{\ell|k+1} \leq (\gamma_u r_u - \rho(1-\lambda^\ell))^2, \; \ell\in\mathbb{N}_1^{N-1}, \\
&\hspace{-0.cm} \phantom{s.t.} \ z_{N|k+1}^\top W_x^{-1} z_{N|k+1} \leq (\gamma_u r_u - \rho(1-\lambda^N))^2, \\
&\hspace{-0.cm} \phantom{s.t.} \ 1 \leq \gamma_u, 
\end{align*}
\end{subequations}
If $\rho$ is such that 
\begin{equation}\label{eq:th1boundrho}
\rho\geq \sqrt{\frac{n(1-\lambda)}{1+\lambda}},
\end{equation}
then, the candidate scaling parameters are expected to not increase over time, i.e., 
\begin{align}
& \E{\gamma_x^+(x_{k},w_k) \, | \, w_k} \leq \gamma_x^\star(x_{k}), \label{eq:condrx}\\
& \E{\gamma_u^+(x_{k},w_k) \, | \, w_k} \leq \gamma_u^\star(x_{k}).\label{eq:condru}
\end{align}
\end{theorem}
The proof can be found in Appendix~\ref{app:rplus}. 

Next, we prove the main result of this paper, that is the stochastic descent property of the cost function, which implies that eventually the closed-loop trajectories will converge to the unconstrained LQR invariant region while the relaxing parameters will converge to one, i.e., the need of constraint relaxation would vanish in expectation.
\begin{theorem}
\label{th:desc}
Consider the system  \eqref{eq:sys_ze} and constraint sets \eqref{eq:XU}, with
$W_x\succ 0$ and $\lambda \in (0,1)$ satisfying \eqref{eq:reachLMI}, and with $r_x$ and $(r_u,W_u)$ computed according to \eqref{eq:Wx} and \eqref{eq:Wu}, respectively. Let $J_N^\star(x_k)$ be the optimal cost for Problem~\eqref{eq:SMPC_mixed}. Then, $J^\star_N(x_k)$ has the following stochastic descent property
\begin{equation*}
\E{J_N^\star(x_{k+1}) \,| \, w_k}-J_N^\star(x_k) \leq \tr{P \Gamma_w}-l(x_k,u^\star_{0|k}),
\end{equation*}
\MF{with $l(x_k,u^\star_{0|k}) = \|x_k\|_Q^2 + \|u^\star_{0|k}\|_R^2$.} Moreover, if there exists a parameter $\mu \in (0,1)$ such that the following conditions hold
\begin{equation}\label{eq:LMIXa}
\frac{W_x^{-1}}{r_{xu}^2} \preceq \frac{Q - \mu P}{\tr{P\Gamma_w}},
\qquad \frac{W_x^{-1}}{r_{xu}^2} \prec \frac{\mu P}{\beta},
\end{equation}
with 
\begin{equation}\label{eq:beta}
\beta = \min\left\{b \in \R^+ \bigg| \ \frac{P}{b} \preceq \frac{Q - \mu P}{\tr{P \Gamma_w}} \right\},
\end{equation} 
then 
\begin{equation}\label{eq:convergence}
\E{J_N^\star(x_{k+i})| \{w_{k+j}\}_{j=0}^{i-1} } \leq (1-\mu)^i J_N^\star(x_k) + \beta/\mu
\end{equation}
for all $i \in \N$, from which it follows
\begin{subequations}\label{eq:limitr}
\begin{align}
& \lim_{i \rightarrow +\infty} \E{x_{k+i} \,|\, \{w_{k+j}\}_{j=0}^{i-1}} \in \Ec_{W_x} (r_{xu}),\\
& \lim_{i \rightarrow +\infty} \E{\g^\star_x(x_{k+i})\,|\, \{w_{k+j}\}_{j=0}^{i-1}}= 1,\\
&\lim_{i \rightarrow +\infty} \E{\g^\star_u(x_{k+i})\,|\, \{w_{k+j}\}_{j=0}^{i-1}} = 1.
\end{align}
\end{subequations}
\end{theorem}
Theorem~\ref{th:desc} proves that the optimal cost of the MS-MPC problem in \eqref{eq:SMPC_mixed} is monotonically decreasing in expectation at every $k$, whenever the state $x_k$ is outside a neighborhood of the origin. Note that, since the cost is composed by a performance and a relaxation term, the decreasing of any of the two terms cannot be ensured at each $k$, as shown in the proof. Nonetheless, from the assessed stochastic descent property, it can be inferred the convergence in expectation of the system state to a probabilistic invariant set and of the relaxation parameter to $1$, meaning that the need of the constraints relaxation will eventually vanish.

\section{Practical issues and probabilistic bounds}
\label{sec:add_res}
This section focuses on two practical aspects related to the application of the proposed MS-MPC approach. First, we show how it is possible to obtain a tighter, a-posteriori probability bound $\rho$ in spite of relaxing the constraint sets, starting from the explicit relation between the relaxation parameters $\gamma_x,\,\gamma_u$ and the violation bound $\varepsilon$. Then, we discuss how to define the bound on the initial input $v_{0|k}$. Indeed, being the first input deterministic, every chance constraint on it would be meaningless in practice. Hence, we outline three different arbitrary strategies that can be implemented and we discuss the corresponding effect on the optimal solution.

\subsection{State-dependent probability bounds}
\label{sec:pr_state}
In Problem~\eqref{eq:SMPC_mixed}, the constraint bounds are softened introducing $\gamma_{x}$ and $\gamma_{u}$ as optimization variables. Alternatively, one could envision to consider parameter $\rho$ as the free optimization variable while maintaining $\gamma_{x} = \gamma_u =1$. This would allow to introduce relaxations on the violation probability bound $\rho$ rather than the constraints relaxation $\gamma_x$, $\gamma_u$. Moreover, since both $\gamma_x,\,\gamma_u$ and $\rho$ appear linearly in Problem~\eqref{eq:SMPC_mixed}, all of them could have been considered as optimization parameters without any relevant complexity increase. In this section, we use the relation between the parameters $\gamma_x$ and $\gamma_u$ with $\varepsilon$ to have tighter, \textit{a-posteriori} estimations of the violation bound $\rho$, starting from the optimal solution to Problem~\eqref{eq:SMPC_mixed}, i.e., $\mathbf{z}^\star_{k},\,\mathbf{v}^\star_{k}$, $\gamma_x^\star$, and $\gamma_u^\star$.

Consider the state constraints first. For each element of the optimal sequence  $z^\star_{\ell|k}$ there exist two possibilities, i.e., either $z_{\ell|k}^\top W_x^{-1} z_{\ell|k} < r_x^2$ or not. In the first case, i.e., $z_{\ell|k}^\top W_x^{-1} z_{\ell|k} < r_x^2$, the nominal state satisfies $z_{\ell|k} \in \Ec_{W_x}\!(r_{x})$ and consequently the probability for $x_{\ell|k} \in \Ec_{W_x}\!(r_{x})$ can be computed by using the Chebychev bound. Hence, it exists a value $\rho_{\ell,x} > 0$ such that 
\begin{align}\label{eq:bound_probability}
\sqrt{z_{\ell|k}^\top W_x^{-1} z_{\ell|k}} = r_{x} - \rho_{\ell,x} (1-\lambda^\ell), 
\end{align}
which implies that $x_{\ell|k} \in \Ec_{W_x}\!(r_{x}) \subseteq \mathbb{X}$ with probability $1-\varepsilon(\rho_{\ell,x})$. Therefore, $1 - \varepsilon(\rho_{\ell,x})$ with 
\begin{equation}\label{eq:barrxk}
\rho_{\ell,x} = \frac{r_{x} - \sqrt{z_{\ell|k}^\top W_x^{-1} z_{\ell|k}}}{1-\lambda^\ell},
\end{equation}
is the maximum lower bound on the probability for $x_{\ell|k}~\in~\Ec_{W_x}\!(r_{x})$, given the specific solution of the optimal problem. 

In the second case, i.e., if $z_{\ell|k}^\top W_x^{-1} z_{\ell|k} \geq r_x^2$, we have that $z_{\ell|k}\not\in\Ec_{W_x}\!(r_{x})$ or it is in its boundary. Therefore, the Chebychev bound cannot be applied to compute a guaranteed probability of constraints satisfaction on $x_k$. In this case, there exists no positive values of $\rho_{\ell,x}$ such that \eqref{eq:bound_probability} holds, since neither the nominal value $z_{\ell|k}$ strictly satisfies the constraints.

Analogous considerations hold for the input constraints. If $v_{\ell|k}^\top W_u^{-1} v_{\ell|k} < r_u^2$ is satisfied, then a guaranteed lower bound on the probability for $u_{\ell|k}~\in~\Ec_{W_u}(r_{u}) \subseteq \U$ is given by $1 - \varepsilon(\rho_{\ell,u})$ with 
\begin{equation}\label{eq:barruk}
\rho_{\ell,u} = \frac{r_{u} - \sqrt{v_{\ell|k}^\top W_u^{-1} v_{\ell|k}}}{1-\lambda^\ell}.
\end{equation}
Finally, if $z_{N|k}^\top W_x^{-1} z_{N|k} < r_{xu}^2 = (\min \{r_x,r_u\})^2$, then $x_{N|k} \in \Ec_{W_x}(r_{xu})$ with  probability at least $1 - \varepsilon(\rho_{xu})$, where
\begin{equation}\label{eq:barrxN}
\rho_{xu} = \frac{r_{xu} - \sqrt{z_{N|k}^\top W_x^{-1} z_{N|k}}}{1-\lambda^N}.
\end{equation}

\subsection{Initial input bound}
\label{sec:inputbounds}
In this section, we discuss three different strategies that can be employed to define the bound on the initial input $v_{0|k}$, each one providing a specific feature to the optimization problem, as detailed in the follows: a) no bounds; b) hard bounds; and c) soft bounds. 

\textit{\underline{Case A - no bounds}}: In case \textit{no bounds} are enforced on $v_{0|k}$, this strategy would lead to more aggressive initial inputs that will drive the nominal states along the prediction to remain inside the region where chance constraints are guaranteed. As a consequence, this approach implies larger input violation at the first step and lower violation occurrences along the (predicted) trajectory. 

\textit{\underline{Case B - hard bounds}}: The second possibility is to enforce \textit{hard bounds} on the first input such that 
$v_{0|k} \in \mathbb{U}$, by simply imposing $H_u v_{0|k} \leq h_u$. This choice would automatically exclude any input constraint violation at the first step, being $v_{0|k}= u_k$ the input to be applied in the MPC strategy control. On the other hand, we could still have input violations along the nominal predicted trajectory for $\ell\in\mathbb{N}_1^{N-1}$. Moreover, this would imply that the recursive feasibility and the estimated tightening bounds are no more valid, since the input violation, admitted within the prediction, is not allowed for the MPC application.  {In other words, the problem solved at time $k+1$ has hard constraints on the first nominal input $v_{0|k+1}$ that were not imposed on $v_{1|k}$.} Hence, the prediction obtained at time $k$ might not be admissible at time $k+1$. This mismatch between prediction and realization appears evidently in the predicted and realized trajectories, depicted in Figure \ref{fig:caseA} and Figure \ref{fig:alltraj}, respectively.

\textit{\underline{Case C - soft bounds}}: The third case, that can be seen as a compromise between the two previous approaches, consists in using the relaxation variable $\gamma_{u}(x_k)$ to \textit{soften} also the bounds on the initial input, i.e., $v_{0|k} \in \gamma_{u}(x_k)\mathbb{U}$. This method implies the possibility of relaxation of the initial input bound, provided that there is no feasible solution for the original constraint, i.e., $\gamma_{u}(x_k)=1$, and the \textit{soft bound} can be obtained by adding to Problem \eqref{eq:SMPC_mixed} the following constraint
\begin{equation}\label{eq:softconstr}
H_u v_{0|k} \leq \gamma_u(x_k)h_u.
\end{equation}

\section{Numerical simulations}
\label{sec:num_sim}
To illustrate the performance of the proposed MS-SMPC, we consider a simple double integrator system 
$$x_{k+1}=Ax_k+Bu_k+w_k,\quad A = \begin{bmatrix}
1 && 1\\
0 && 1
\end{bmatrix}, \quad 
B = \begin{bmatrix}
0.5\\
1
\end{bmatrix},
$$
also used as benchmark example in \cite{hewing2018stochastic,hewing2020recursively,arcari2023stochastic}. Here, the local feedback gain $K = \left[0.2068 \ \, 0.6756\right]$ is the infinite horizon, (discrete) LQR solution with $Q = \mathbb{I}_2$, $R = 10$ such that \MF{$P=\begin{bmatrix}
    3.2664 & 3.2016; & 3.2016 & 9.3569\end{bmatrix}$.}
The covariance matrix $\Gamma_w$ of the i.i.d. Gaussian process $w_k$, with null mean, and the matrix $W_x$, which satisfies the conditions \eqref{eq:reachLMI_1} and \eqref{eq:reachLMI_2} with $\lambda = 0.7503$, are given by 
\begin{equation*}
\Gamma_w = \left[\begin{array}{cc}
0.1 & 0.05\\
0.05 & 0.1
\end{array}\right], \quad W_x = \left[\begin{array}{cc}
10.9264  & -3.7386\\
-3.7386  &  3.8143
\end{array}\right].
\end{equation*}
The polytopic state and input constraint sets are
\begin{align*}
& \X = \{x \in \R^2: H_x x \leq h_x\} = \{x \in \R^2: \|x\|_\infty \leq 40\}\\
& \U = \{u \in \R: H_u u \leq h_u\} = \{u \in \R: |u| \leq 10\}
\end{align*}
giving $r_{x} = 12.1010$ as the maximal value such that $\Ec_{W_x}\!(r_{x}) \subseteq \X$, and $W_u = 0.2237$ with $r_{u} = 21.1448$ the maximal value such that $\Ec_{W_u}\!(r_u) \subseteq \U$ and \eqref{eq:WuWx} is satisfied. The selected violation probability level is $\varepsilon~=0.1$, leading to $\rho = 2.146$ that satisfies \eqref{eq:th1boundrho}, while $\eta = 10^5$ is used in the cost of (\ref{eq:SMPC_mixed}). Conditions (\ref{eq:LMIXa}) for convergence in expectation, i.e., (\ref{eq:limitr}), hold with $\mu =  0.0464$ and $\beta = 33.7956$.

{For comparison purposes, we evaluate our approach against a standard dual-mode open-loop strategy, in the spirit of \cite{farina2013probabilistic,hewing2018stochastic}. In particular, we consider the following deterministic problem
\begin{equation}
\label{eq:SMPC_init}
\begin{aligned}
\min_{\textbf{z}_k, \textbf{v}_k} &\sum_{\ell = 0}^{N-1} \left(\|z_{\ell|k}\|_Q^2+ \|v_{\ell|k}\|_R^2\right) + \|z_{N|k}\|_P^2\\
\text{s.t.} \ \ & z_{\ell+1|k} = A z_{\ell|k} + B v_{\ell|k},\\
& z_{\ell|k}\in\mathbb{Z}_\ell \subseteq \X \ominus \mathcal{R}_{\ell}^x, \quad \ell\in\mathbb{N}_1^{N-1}\\
& v_{\ell|k}\in\mathbb{V}_\ell \subseteq \U \ominus \mathcal{R}_{\ell}^u, \quad \ell\in\mathbb{N}_1^{N-1}\\
& z_{N|k} \in \mathbb{Z}_N \subseteq \XN \ominus \mathcal{R}_{\infty}^x,
\end{aligned}
\end{equation}
with probabilistic sets defined as \eqref{eq:PRS_ell} and \eqref{eq:PRS_ell_u}, and with $\mathcal{R}_\infty^x = \Ec_{W_x}\!\big(\rho\big)$, and considering a standard dual-mode initialization approach, i.e., $z_{0|k}=x_k$ if \eqref{eq:SMPC_init} is feasible for $x_k$ or  $z_{0|k}=z_{1|k-1}$ otherwise. Note that for this problem, referred to as initial-state SMPC (IS-SMPC) in the follows, the recursive feasibility conditions is ensured only in if it is feasible at the initial time.}

\subsection{Comparison of input bounds strategies}
In the first set of simulations, to better highlight the effect of the proposed relaxation of {MS-SMPC} with respect to {IS-SMPC}, we select  as initial state the point $x_0 = (-40, \, 40)$, which is on the boundary of the constraint set $\X$. 
\begin{figure}[!ht]
\centering
\subfigure[case A]{
\includegraphics[width=\columnwidth]{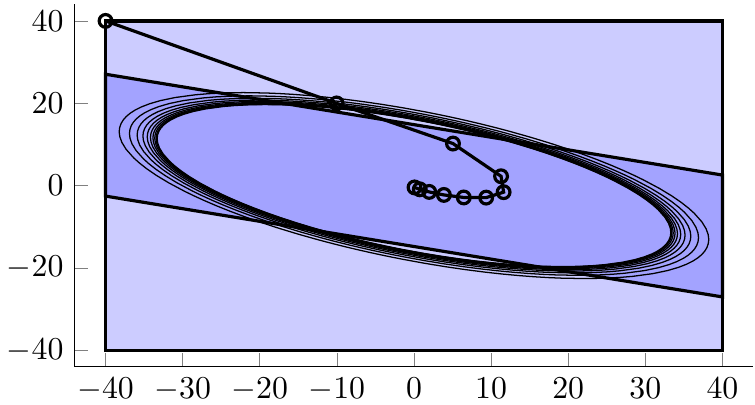}\label{fig:caseA0}
}
\hfil
\subfigure[case B]{
\includegraphics[width=\columnwidth]{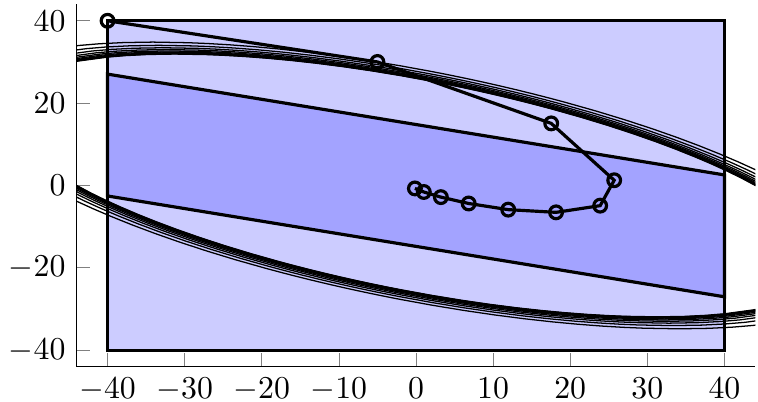}\label{fig:caseB0}
}
\subfigure[case C]{
\includegraphics[width=\columnwidth]{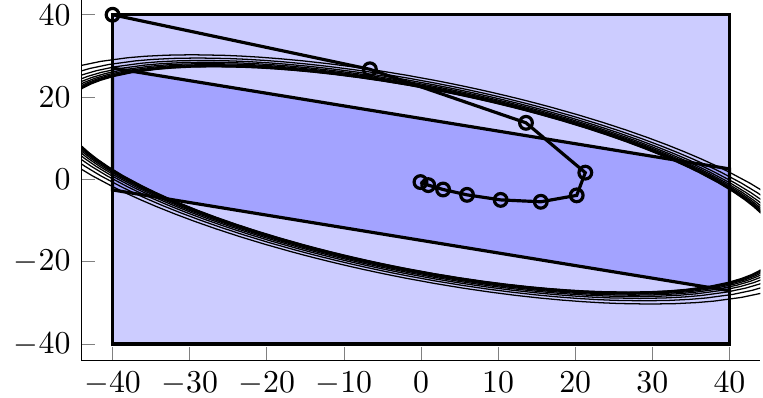}\label{fig:caseC0}
}
\caption{Predicted state trajectories and sets $\mathbb{Z}_\ell$, starting from $x_0=(-40,\,40)$ with $\ell\in\mathbb{N}_0^{10}$, enforcing on $v_{0|k}$: no bounds (top), hard bounds (middle), and soft bounds (bottom).}
\label{fig:caseA}
\end{figure}
This selected initial condition leads to infeasibility of {IS-SMPC} and therefore the problem cannot be tackled by classical open-loop strategies. On the other hand, the proposed {MS-SMPC} is able to overcome the feasibility issue through the relaxation of the constraints by employing values of $\gamma_{x}$ and $\gamma_{u}$, greater than one, in \eqref{eq:SMPC_mixed}. 

Now, we proceed analyzing the effects of the different assumptions on the initial input strategies described in Section \ref{sec:inputbounds}, when applied to the aforementioned system, namely by enforcing no bounds (case A), hard bounds (case B), or soft bounds (case C) on $v_{0|k}$ and solving the related optimization problem \eqref{eq:SMPC_mixed} with $N = 10$. 

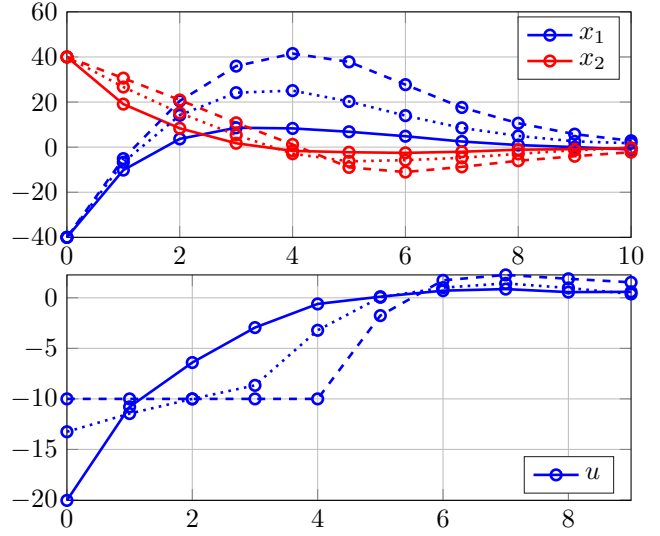
\begin{figure}[!ht]
\begin{center}
%
%
\begin{tikzpicture}

\begin{axis}[%
width=7.5cm,
height=3cm,
at={(0.0cm,3.5cm)},
scale only axis,
xmin=0,
xmax=10,
ymin=-40,
ymax=60,
axis background/.style={fill=white},
xmajorgrids,
ymajorgrids,
legend style={legend cell align=left, align=left, draw=white!15!black}
]
\addplot [color=blue, solid, line width=1pt, mark=o, mark options={solid, blue}]
  table[row sep=crcr]{%
0	-40\\
1	-10.1796913800563\\
2	3.66018700947488\\
3	8.60835206260553\\
4	8.31633238372446\\
5	6.82813694083673\\
6	4.84839464203473\\
7	2.5090808960426\\
8	0.929558406056311\\
9	0.0184581349622375\\
10	-1.03039460380614\\
};
\addlegendentry{$x_1$}

\addplot [color=red, solid, line width=1pt, mark=o, mark options={solid, red}]
  table[row sep=crcr]{%
0	40\\
1	19.0862207624914\\
2	8.35818146501316\\
3	1.72328041405526\\
4	-1.65081642433511\\
5	-2.25478089492298\\
6	-2.53163235958191\\
7	-2.04763209482453\\
8	-1.12098620753253\\
9	-0.882996117259601\\
10	-0.397975803014893\\
};
\addlegendentry{$x_2$}

\addplot [color=blue, dashed, line width=1pt, mark=o, mark options={solid, blue}, forget plot]
  table[row sep=crcr]{%
0	-40\\
1	-5.08626630240261\\
2	20.4849209192965\\
3	35.9421217595103\\
4	41.4924823273714\\
5	37.8277750287453\\
6	27.7211175668035\\
7	17.5893884384798\\
8	10.6086814755428\\
9	5.76619086672899\\
10	2.82403977311957\\
};
\addplot [color=red, dashed, line width=1pt, mark=o, mark options={solid, red}, forget plot]
  table[row sep=crcr]{%
0	40\\
1	30.6032698655761\\
2	20.9091238068917\\
3	10.7920880472217\\
4	1.20559579092008\\
5	-8.97730605552702\\
6	-11.0481259370345\\
7	-8.73927163000825\\
8	-6.03591753577041\\
9	-4.03782237492779\\
10	-2.21307177532432\\
};
\addplot [color=blue, dotted, line width=1pt, mark=o, mark options={solid, blue}, forget plot]
  table[row sep=crcr]{%
0	-40\\
1	-6.71087970901678\\
2	14.134855463645\\
3	24.1549355655325\\
4	25.0747291627431\\
5	20.248719801953\\
6	13.9811066381559\\
7	8.53311393979558\\
8	4.97849382651633\\
9	2.57243999558841\\
10	1.73203973343294\\
};
\addplot [color=red, dotted, line width=1pt, mark=o, mark options={solid, red}, forget plot]
  table[row sep=crcr]{%
0	40\\
1	26.6652409784717\\
2	15.1423963937701\\
3	5.4258233437534\\
4	-2.92127962826978\\
5	-6.27406840654988\\
6	-5.76926520626399\\
7	-4.72383963759547\\
8	-2.96840645247643\\
9	-1.35511722540379\\
10	-0.389600049526245\\
};
\end{axis}

\begin{axis}[%
width=7.5cm,
height=3cm,
at={(0cm,0cm)},
scale only axis,
xmin=0,
xmax=9,
ymin=-20.0460977152023,
ymax=2.26619300917832,
axis background/.style={fill=white},
xmajorgrids,
ymajorgrids,
legend pos=south east,
legend style={legend cell align=left, align=left, draw=white!15!black}
]
\addplot [color=blue, solid, line width=1pt, mark=o, mark options={solid, blue}]
  table[row sep=crcr]{%
0	-20.0460977152023\\
1	-10.7891580950678\\
2	-6.40432177129469\\
3	-2.94470107973098\\
4	-0.604795979458538\\
5	0.111066238622278\\
6	0.707630361594869\\
7	0.864506731832384\\
8	0.565114674717901\\
9	0.59276777800796\\
};
\addlegendentry{$u$}

\addplot [color=blue, dashed, line width=1pt, mark=o, mark options={solid, blue}, forget plot]
  table[row sep=crcr]{%
0	-9.9999999542667\\
1	-9.99999992238295\\
2	-9.99999997072594\\
3	-9.99999998048203\\
4	-9.9999999372548\\
5	-1.7591054757517\\
6	1.73056532076997\\
7	2.26619300917832\\
8	1.8838832025205\\
9	1.53525423349608\\
};
\addplot [color=blue, dotted, line width=1pt, mark=o, mark options={solid, blue}, forget plot]
  table[row sep=crcr]{%
0	-13.2546999979462\\
1	-11.4543439899055\\
2	-9.99999995807063\\
3	-8.66168359790677\\
4	-3.21259351597071\\
5	0.05076179970407\\
6	1.00600757783647\\
7	1.42652968811479\\
8	0.975745059911104\\
9	0.383458676899916\\
};
\end{axis}
\end{tikzpicture}%
\end{center}
\caption{State and input trajectories for SMPC enforcing no bounds (solid lines), hard bounds (dotted lines), and soft bounds (dashed lines) on the first input $v_{0|k}$.}
\label{fig:allstatein} 
\end{figure}

Figure \ref{fig:caseA} depicts the predicted trajectories given by $z_{\ell|k}$ and the sets $\mathbb{Z}_\ell=\Ec_{W_x}\!(\gamma_x r_x - \rho(1-\lambda^\ell))$ for $\ell \in \N_1^{N}$ for the three cases. In case of \textit{no constraint}, the {MS-SMPC} generates a predicted trajectory that reaches quickly the set $\X_K$ where both state and control input are satisfied (see Figure~\ref{fig:caseA0}), and the value $\gamma_{x}$ and $\gamma_{u}$ are close to one. The drawback is a non negligible violation on the input constraint, as shown in Figure~\ref{fig:allstatein}.

On the other hand, the  {MS-SMPC} with \textit{hard bounds} on the first predicted input (case B) leads to a slower convergent trajectory, as depicted in Figure~\ref{fig:caseB0}, and necessitates of higher relaxations on the state constraints, as witnessed by the larger sequence of ellipsoidal sets. However, this approach ensures no input constraints violation (see Figure~\ref{fig:allstatein}).

Finally, the \textit{soft input strategy} (case C), which is a trade-off between the other two, leads to a state trajectory that is faster than the one obtained with hard constraints, but with more input violations occurrences, as illustrated in Figure~\ref{fig:caseC0}. Moreover, as depicted in Figure~\ref{fig:allstatein}, the input constraints are violated in the first three steps but with lower values, with respect to the strategy of no bounds, leading to slower prediction and then higher state constraints violations probabilities.

\begin{figure}[!ht]
\begin{center}
\includegraphics[width=\columnwidth]{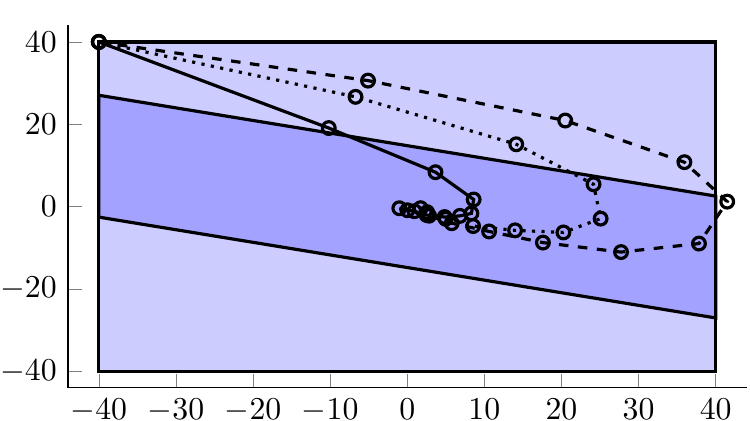}
\end{center}
\caption{Comparison among the realized state trajectories starting from $x_0=(-40,\,40)$ for $k \in [0,10]$ enforcing on $v_{0|k}$: a) no bounds (solid line); b) hard bounds (dotted line); and, c) soft bounds (dashed lined).}
\label{fig:alltraj} 
\end{figure}

The considerations above are corroborated also by the results depicted in Figure~\ref{fig:alltraj}, representing the realized trajectories obtained along a simulation horizon of 10 steps for the three strategies. It can be noticed that the tighter the constraints on the initial input are, the slower the trajectory convergence rate is. 

\begin{figure}[!ht]
\begin{center}
\includegraphics[width=\columnwidth]{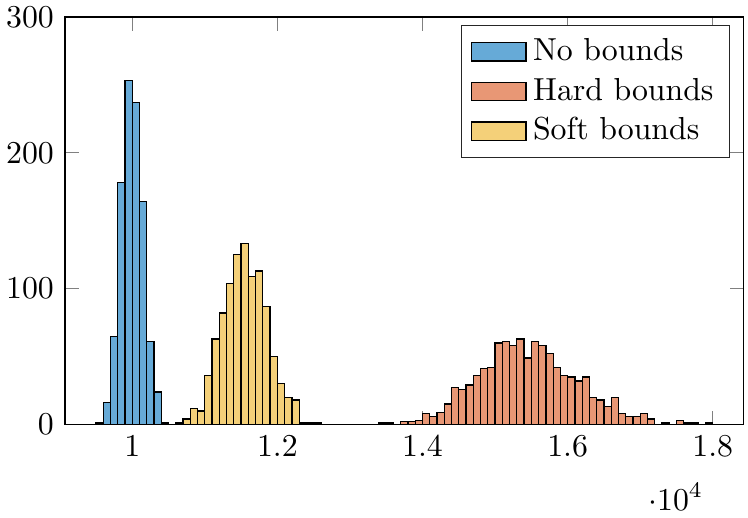}
\caption{Cost comparison over 1000 simulations with $N=10$ for $x_0 = (-40,\, 40)$ of the three input bounds strategies. Mean cost: no bounds $J_{MPC} = 9999$; hard bounds $J_{MPC}=15460$; soft bounds $J_{MPC}=11552$. \label{fig:hists4040}}
\end{center}
\end{figure}

Last, the three strategies have been compared also in terms of the following performance index
\begin{equation}\label{eq:costMPC}
J_{MPC} \doteq\,\, \sum_{k = 0}^{N-1} \left(\|x_k\|_Q^2+ \|u_k\|_R^2\right)
\end{equation}
over 1000 trajectories with $N=10$ and $x_0=(-40,\,40)$. The histograms of the resulting cost values are depicted in Figure~\ref{fig:hists4040}. It can be noticed that, the tighter the bound on the initial input $v_{0|k}$ is, the worst the control performances are, as expected.

\subsection{Probability bounds comparison}
As illustrated in Section~\ref{sec:pr_state}, the constraints relaxation obtained by introducing $\gamma_{x}$ and $\gamma_{u}$ as free variables in the optimization problem \eqref{eq:SMPC_mixed} can also be interpreted as a relaxation in the probability levels of the chance constraints. In particular, by determining $\rho_{\ell,x}$ as in \eqref{eq:barrxk} and $\rho_{\ell,u}$ as in \eqref{eq:barruk} for a solution $z_{\ell|k}$ and $v_{\ell|k}$ with $\ell = \mathbb{N}_1^N$, the probability bounds on the chance constraints satisfaction can be obtained for the specific measured state $x_k$. Table~\ref{tab:1} provides the probability bound values 
\begin{equation}\label{eq:rho}
p_\ell(j) = \big(p_{\ell,x}(j), \, p_{\ell,u}(j)\big) = \Big( \chi_n^2 \big(\rho_{\ell,x}^2(j) \big) , \, \chi_n^2 \big(\rho_{\ell,u}^2(j)\big) \Big)
\end{equation}
with $\rho_{\ell,x}(j)$ and $\rho_{\ell,u}(j)$ computed as in \eqref{eq:barrxk}--\eqref{eq:barruk} for each input bound strategy, i.e., $j = A, B, C$, and for all $\ell = \mathbb{N}_1^{10}$, starting from the optimal prediction pair, $z_{\ell|k}$ and $v_{\ell|k}$, computed at $k = 0$. Moreover, $N_{sim} = 1000$ simulations have been run for the different strategies and the occurrence of violations of the constraints $x_k \in \Ec_{W_x}\!(r_{x})$ and $u_k \in \Ec_{W_u}\!(r_{u})$ are registered to obtain the relative frequency of constraints satisfaction given by
\begin{align}
f_\ell(j) & = \big(f_{\ell,x}(j), \, f_{\ell,u}(j)\big) \nonumber\\
& = \left( \frac{N_{sim} - n_{\ell,x}(j)}{N_{sim}}, \ \frac{N_{sim} - n_{\ell,u}(j)}{N_{sim}} \right) \label{eq:f}
\end{align}
where $n_{\ell,x}(j)$ and $n_{\ell,u}(j)$ are the number of violations at time $\ell$ of the state and input constraints, respectively, with $j = A,B,C$. The values of $f_\ell(j)$ are reported in the right part of Table~\ref{tab:1}. It is worth noting that bounds on the future constraints satisfaction, based on the optimal nominal states and inputs obtained at time $k = 0$, are reasonably accurate guesses. 

Moreover, we can note how the constraints on the first input $v_{0|k}$ affect the violation probabilities of the MS-SMPC. Larger freedom on the first input selection leads to a more aggressive MPC action for driving the state towards the feasibility region with no relaxation, hence to bigger violations on the input constraints (see Figure~\ref{fig:allstatein}) but less frequent violation on the future state constraints, which is also reasonable.

\begin{table}[!ht]
\centering
\setlength\tabcolsep{3pt} 
\begin{tabular}{|c|ccc|ccc|}
 \hline \hline
$\ell$ &  $p_\ell(A)$ & $p_\ell(B)$ & $p_\ell(C)$ & $f_\ell(A)$ &  $f_\ell(B)$ & $f_\ell(C)$\\
\hline
1  & (1,0.90)	& (0,0) 		& (0,0) 	& (1,0)	  	& (0,0)		& (0,0)     \\
2  & (1,1) 		& (0, 0) 		& (0,0) 	& (1,1)    	& (0,0)     & (0,0)   	\\
3  & (1,1) 		& (0.99,1)	 	& (1,1) 	& (1,1)    	& (0,0)  	& (0.97,1)  \\
4  & (1,1) 		& (1, 1) 		& (1,1) 	& (1,1)    	& (0,0.89)  & (1,1)     \\
5  & (1,1) 		& (1,1)		  	& (1,1)		& (1,1)    	& (0.97,1)  & (1,1)     \\
6  & (1,1) 		& (1,1)   		& (1,1)		& (1,1)    	& (1,1)     & (1,1)     \\
7  & (1,1)		& (1,1)   		& (1,1)		& (1,1)    	& (1,1)     & (1,1)     \\
8  & (1,1)		& (1,1)   		& (1,1)		& (1,1)    	& (1,1)     & (1,1)     \\
9  & (1,1)		& (1,1)   		& (1,1)		& (1,1)    	& (1,1)     & (1,1)     \\
10 & (1,1)    	& (1,1)       	& (1,1)     & (1,1)		& (1,1)   	& (1,1)\\
\hline \hline
\end{tabular}
\caption{Probability of state and input constraints satisfaction within the predicted bounds $p_\ell(j)$ defined in \eqref{eq:rho} with $\rho_{\ell,x}(j)$ and 
$\rho_{\ell,u}(j)$ as in \eqref{eq:barrxk} and \eqref{eq:barruk} for $j = A, B, C$, and  relative frequencies $f_\ell(j)$ defined in \eqref{eq:f} of both state and input constraints satisfaction along the trajectories for MS-SMPC over $N_{sim} = 1000$ tests each, for $\ell \in \N_1^{10}$ and $x_0 = (-40, \, 40)$. \label{tab:1}
}
\end{table}

\begin{figure}
\begin{center}
\includegraphics[width=\columnwidth]{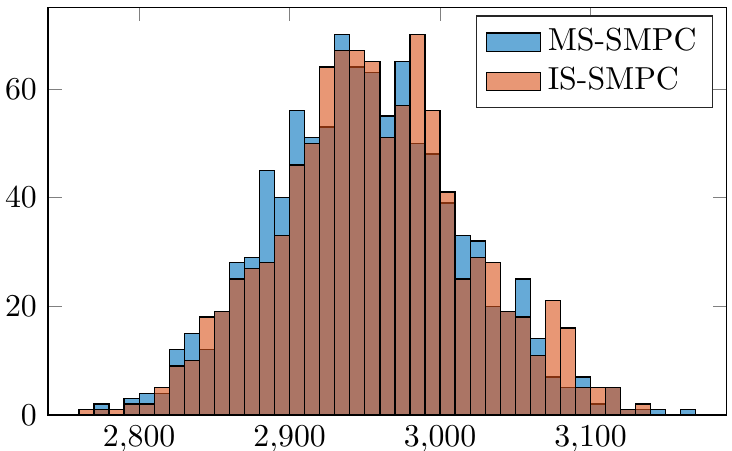}
\caption{Cost comparison over 1000 simulations with $N=10$ for $x(0) = (-30, 0)$. Mean cost: for the new  {MS-SMPC} (with no initial bound) $J_{MPC}=2951$, for the  {IS-SMPC} $J_{MPC}=2956$, with an almost unitary ratio. \label{fig:3000}}
\end{center}
\end{figure}

\subsection{Comparison of  {IS-SMPC} and  {MS-SMPC}}
In this last set of simulations, we consider a case where both IS-SMPC and MC-SMPC are feasible, and we compare the  {MS-SMPC} scheme (with no bounds on in the initial input) with the  {IS-SMPC} approach in terms of the performance index $J_{MPC}$ defined in \eqref{eq:costMPC}. For the first case study, we consider as initial condition $x(0) = (-30,0)$, which is inside the feasibility region and far from its boundaries. For each SMPC scheme, $N_{sim}=1000$ simulations have been run. In Figure~\ref{fig:3000} we can observe that both control schemes lead, in practice, to the same unconstrained optimization problem and, consequently, to the same solutions. Indeed, the histograms of the cost function realized along the generated trajectories and depicted in Figure~\ref{fig:3000} are almost overlapped, and the mean value of the performance index are almost identical.

\begin{figure}
\begin{center}
\includegraphics[width=\columnwidth]{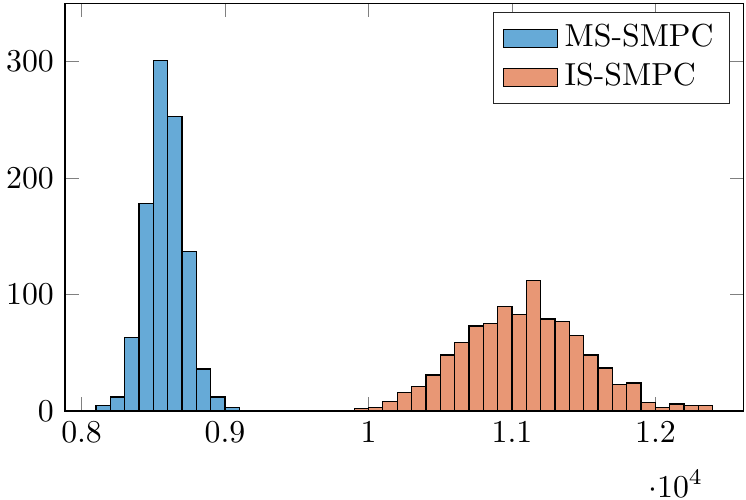}
\caption{Cost comparison over 1000 simulations of 10 steps for $x(0) = (-40, 37)$, close to the  {IS-SMPC} feasibility bounds. Mean cost: for the  {MS-SMPC} (with no initial bound) $8584$, for the IF-SMPC $11085$, with ratio of around $0.77$. \label{fig:4037}}
\end{center}
\end{figure}

On the other hand, for initial conditions closer to the boundary of the feasibility region, we can observe quite distinct behaviors of the two SMPC schemes. In particular, running $N_{sim}=1000$ simulations starting from $x(0) = (-40,37)$, we can notice a significantly different effect of the constraints on the realized trajectories, as shown in Figure~\ref{fig:4037}. Indeed, in this case study, the histograms representing the cost $J_{MPC}$ for the obtained realizations show that the performance for the  {MS-SMPC} substantially outperforms the one granted by the  {IS-SMPC}, being the mean of the former less than $80\%$ of the latter.

\section{Conclusion and future works}
\label{sec:concl}

Motivated by the difficulty of guaranteeing recursive feasibility of SMPC problems in the presence of unbounded stochastic  uncertainties, especially when one wants to exploit the knowledge of the measured state in implementing the
feedback strategy, in this paper we
introduced a novel definition of recursive feasibility. Namely, we allow for a relaxation of the constraints so that the ensuing optimization problem remains feasible, but we require that such relaxation is minimal, and more importantly, it has certain probabilistic properties. In parallel, we developed a novel SMPC approach guaranteeing these new properties. We prove, also by numerical simulations, that this new approach allows to tackle problems not affordable by classical open-loop strategies, and that it outperforms these latter when both are feasible.

In our next studies, we aim to extend the proposed philosophy to other SMPC frameworks. For instance, we will investigate how the introduced concepts may applied to problems involving {mission-wide constraints}, and more specifically \textit{mission-wide probability of safety} (see e.g., \cite{Lew2019ChanceConstrainedOA, Gross2021mission, raghuraman2022long}). Moreover, we will investigate the extension of the proposed approach beyond the ellipsoidal approximation, where the concept of relaxation and the associated probabilistic properties will be revised to involve polytopic constraint sets.
\appendix
\section*{Appendix}

\section{Proof of Proposition \ref{prop:xellk}}
\label{app:prop1}
First, we notice that, by definition, $z_{\ell|k} \in \Ec_{W_x}\!(\alpha_\ell)$ with $\alpha_\ell = \Big(z_{\ell|k}^\top W_x^{-1}z_{\ell|k}\Big)^{1/2}$. Then, combining  \eqref{eq:PRS_ell} with the relation $\mathbb{Z}_\ell\subseteq \X\ominus \mathcal{R}_\ell^x$, one obtains
\begin{equation}\label{eq:prop1b}
     \Pr{x_{\ell|k} \in \Ec_{W_x}\Big(\alpha_\ell+\rho(1-\lambda^\ell)\Big)}\geq 1-\varepsilon.
 \end{equation}
Then, to ensure \eqref{eq:Prop1a}, being $\Ec_{W_x}\!(r_x)  \subseteq\mathbb{X}$ by construction, it is sufficient to prove the following set inclusion
\begin{equation}
   \Ec_{W_x}\Big( \alpha_\ell+\rho(1-\lambda^\ell)\Big)  \subseteq \Ec_{W_x}\!(r_x),
\end{equation}
which holds if and only if $\alpha_\ell+\rho(1-\lambda^\ell) \leq r_x$ or, equivalently, if $z_{\ell|k}^\top W_x^{-1}z_{\ell|k}\leq \Big(r_x-\rho(1-\lambda^\ell)\Big)^2$, that is if \eqref{eq:prop1} is satisfied.

\section{Proof of Proposition \ref{prop:uellk}}
\label{app:prop2}
Let us consider the constraints on the nominal input $v_{\ell|k}$. By construction, see \eqref{eq:WuWx}, if $e_{\ell|k} \in \Ec_{W_x}\!\left(\rho(1-\lambda^\ell)\right)$ then $K e_{\ell|k} \in \mathcal{R}_\ell^u = \Ec_{W_u}(\rho (1-\lambda^\ell))$. Since $e_{\ell|k} \in \Ec_{W_x}\!\left(\rho(1-\lambda^\ell)\right)$ holds with probability $1-\varepsilon$ by definition, being the set  $\Ec_{W_x}\!\left(\rho(1-\lambda^\ell)\right)$ the PRS at time $\ell$ for \eqref{eq:sys_e}, it follows that 
\begin{equation}\label{eq:prop2a}
v_{\ell|k} \in \Ec_{W_u}(r_u - \rho (1-\lambda^\ell))
\end{equation}
implies $u_{\ell|k} = v_{\ell|k} + K e_{\ell|k} \in \Ec_{W_u}(r_u)$ with the same probability $1-\varepsilon$. Hence, the input chance constraints are satisfied with probability no smaller than $1-\varepsilon$ if 
\begin{equation*}
v_{\ell|k}^\top W_u^{-1} v_{\ell|k}\leq (r_{u}- \rho (1-\lambda^\ell))^2.
\vspace*{-0.5cm}
\end{equation*}

\section{Proof of Lemma \ref{lem:rW}}
\label{app:lemma1}
From \eqref{eq:WuWx}, we have that $x \in \Ec_{W_x}(r)$ implies $Kx \in \Ec_{W_u}(r)$ for all $r \geq 0$. Hence, for $r=r_{xu}$, we have  
$x \in \Ec_{W_x}(r_{xu})$ implies $Kx \in \Ec_{W_u}(r_{xu}).$ Now, let us first consider the case $r_{xu} = r_x \leq r_u$. Then, we have $x \in \Ec_{W_x}(r_{xu}) = \Ec_{W_x}(r_x) \subseteq \X,$ and $u = Kx \in \Ec_{W_u}(r_{xu}) \subseteq \Ec_{W_u}(r_u) \subseteq \U.$ Analogously, for $r_{xu} = r_u < r_x$, then $x \in \Ec_{W_x}(r_{xu}) \subset \Ec_{W_x}(r_x) \subseteq \X,$ and $u = Kx \in \Ec_{W_u}(r_{xu}) = \Ec_{W_u}(r_u) \subseteq \U.$

\section{Proof of Theorem~\ref{th:term}}
\label{app:th_term}
First, we notice that
\begin{equation}\label{eq:th1}
z_{N|k} \in \Ec_{W_x}\!(\alpha_N), \quad \text{with} \quad \alpha_N = \Big(z_{N|k}^\top W_x^{-1} z_{N|k}\Big)^{1/2}.
\end{equation}
Moreover, since $\Ec_{W_x}\!\left(\rho(1-\lambda^N)\right)$ is a probabilistic reachable set from \eqref{eq:reachLMI}, we have that
\begin{equation}
\Pr{e_{N|k} \in \Ec_{W_x}\!\left(\rho(1-\lambda^N)\right)} \geq 1-\varepsilon,
\end{equation}
which implies that
\begin{equation}\label{eq:condN}
\Pr{x_{N|k} \in \Ec_{W_x}\!\left( \alpha_N + \rho (1-\lambda^N) \right)} \geq 1-\varepsilon.
\end{equation}
Hence, similarly to Proposition~\ref{prop:xellk}, see \eqref{eq:prop1}, and relying on the definition of PRS, we have that \eqref{eq:constr_zNk} implies 
\begin{equation}\label{eq:th1b}
\Ec_{W_x}\!\left( \alpha_N + \rho (1-\lambda^N)\right) \subseteq \Ec_{W_x}\!\left( r_{xu}\right) \subseteq \X_N,
\end{equation}
guaranteeing chance constraints satisfaction for $\ell=N$.

Then, from \cite[Proposition 4]{fiacchini2021probabilistic} combining \eqref{eq:Pre} with
\begin{equation}\label{eq:lambda}
A_K^\top W_x^{-1} A_K \preceq \lambda^2 W_x^{-1},
\end{equation}
for all $j \in \N$ we have
\begin{align}
&z_{N+j|k} \in \Ec_{W_x}\!(\alpha_N \lambda^j), \\
&\Pr{e_{N+j|k} \in \Ec_{W_x}\!\left( \rho(1-\lambda^{N+j})\right) } \geq 1-\varepsilon, \label{eq:th1_a}
\end{align}
which implies that 
\begin{equation}\label{eq:condK1}
\Pr{x_{N+j|k} \in \Ec_{W_x}\!\left(\alpha_N \lambda^j + \rho(1-\lambda^{N+j}) \right)} \geq 1-\varepsilon.
\end{equation}
It is now left to prove that \eqref{eq:condK1} implies \eqref{eq:constr_xkx} for all $\ell \geq N$ (for $\g_x=1$). First, consider the case $\alpha_N <\rho \lambda^N$. Hence, we have
\begin{equation*}\label{eq:th1c1}
\Ec_{W_x}\!\left( \alpha_N \lambda^j + \rho(1-\lambda^{N+j})\right) \subseteq \Ec_{W_x}\!\left(\rho \right)
\end{equation*}
and then, from $\rho \leq r_{xu}$ and \eqref{eq:condK1}, we get
\begin{equation*}\label{eq:th1c}
\Ec_{W_x}\!\left( \alpha_N \lambda^j + \rho(1-\lambda^{N+j})\right) \subseteq \Ec_{W_x}\!\left(\rho \right) \subseteq \Ec_{W_x}\!\left( r_{xu}\right) \subseteq \XN,
\end{equation*}
implying constraints satisfaction with $v_{\ell|k} = K z_{\ell|k}$.
Now, considering the case 
\begin{equation}\label{eq:cond_decreas}
\alpha_N \geq \rho\lambda^{N},
\end{equation}
we search for a condition ensuring that \eqref{eq:condN} implies 
\begin{equation}\label{eq:condK}
\Pr{x_{N+j|k} \in \Ec_{W_x}\!\left(\alpha_N + \rho(1-\lambda^N) \right)} \geq 1-\varepsilon, 
\end{equation}
for all $j \in \N$, see \eqref{eq:th1b}.
%
To guarantee the satisfaction of condition \eqref{eq:condK}, from \eqref{eq:condK1},  it is sufficient to have 
\begin{equation}
\alpha_N \lambda^j+ \rho(1-\lambda^{N+j}) \leq \alpha_N + \rho(1 - \lambda^{N}),
\end{equation}
for all $j \in \N$, condition holding if and only if \eqref{eq:cond_decreas} is satisfied.
Moreover, $\alpha_N$ satisfying \eqref{eq:constr_zNk} is equivalent to 
\begin{equation}\label{eq:cond_x11}
\alpha_N \leq r_{xu} - \rho(1-\lambda^N).
\end{equation}
Therefore, combining \eqref{eq:cond_decreas} and \eqref{eq:cond_x11}, we obtain that $\Pr{x_{N+j|k} \in \XN} \geq 1-\varepsilon$ is satisfied for all $j \in \N$ if 
\begin{equation}
\lambda^{N} \rho \leq \alpha_N \leq r_{xu} - \rho(1-\lambda^N), 
\end{equation}
which has an admissible solution if and only if \eqref{eq:constr_zNk} is satisfied with $\rho \leq r_{xu}$, assumed holding. Consequently, to guarantee that both state and input chance constraints are satisfied by the system with control law $u_{\ell|k} = K x_{\ell|k}$ along the whole future trajectory starting at time $N$, it is sufficient to have \eqref{eq:constr_zNk} as terminal constraint for the nominal system, {from the definition \eqref{eq:XN} of $\XN$.}


\section{Proof of Proposition~\ref{prop:inEcWx}}
\label{app:Prop3}
First, we consider the nominal state evolution and the related constraint sets along the horizon $N$. For $z_{0|k}\in\Ec_{W_x}(r_{xu})$, we have that
$$z_{0|k}^\top W_x^{-1}z_{0|k}\leq (\gamma_x r_{xu})^2.$$
Then, for $z_{1|k}$ with $v_{0|k}=Kz_{0|k}$, and considering \eqref{eq:reachLMI_1} we have $z_{1|k}=A_Kz_{0|k}$ and consequently
$$z_{1|k}^\top A_K^\top W_x^{-1}A_K z_{1|k}\leq \lambda^2(\gamma_x r_{xu})^2.$$
Moreover, since $\rho\leq r_{xu}\leq \gamma_x r_{xu}$, then
$$(\lambda \gamma_x r_{xu})^2=(\gamma_x r_{xu}-\gamma_x r_{xu}(1-\lambda))^2\leq (\gamma_x r_{xu}-\rho(1-\lambda))^2.$$
Hence, for $z_{1|k}$, we obtain that
\begin{align*}
& z_{1|k}^\top A_K^\top W_x^{-1} A_K z_{1|k} \leq \lambda^2(\gamma_x r_{xu}-\rho(1-\lambda))^2\\
&=\big(\lambda(\gamma_x r_{xu}-\rho)+\lambda^2\rho\big)^2 \leq \big( \gamma_x r_{xu}-\rho+\lambda^2\rho\big)\\
&=\big( \gamma_x r_{xu}-\rho(1-\lambda^2)\big).
\end{align*}
Proceeding analogously for $\ell>2$, it can be proved that for any $\ell$ we have that, if $x_k\in\Ec_{W_x}(r_{xu})$, we have
$$z_{\ell|k}^\top W_x^{-1}z_{\ell|k}\leq \big(\gamma_x r_{xu}-\rho(1-\lambda^\ell)\big)^2.$$
Analogous considerations lead to prove that for all $\ell\in\N$, if $x_k\in\Ec_{W_x}(r_{xu})$ then
$$v_{\ell|k}^\top W_u^{-1} v_{\ell|k}=z_{\ell|k}^\top K^\top W_u^{-1}Kz_{\ell|k}\leq \big(\gamma_u r_{xu}-\rho(a-\lambda^\ell)\big)^2.$$
Hence, if $\gamma_x^\star=\gamma_u^\star=1$, then Problem~\eqref{eq:SMPC_mixed} is a unconstrained LQR problem, which concludes the proof.

\section{Proof of Theorem \ref{th:rplus}}
\label{app:rplus}
Note first that state $x_{k+1}$ predicted at time $k$ based on $x_k$ and the realization of $w_k$ is a random variable with mean $\E{x_{k+1}| w_k} = A x_k + B v_{0|k}^\star(x_k)$ and covariance $\Gamma_w$, being $x_{k+1} = \E{x_{k+1} | w_k} + w_k$. Hence, also $\z^+$, $\v^+$, $\gamma_x^+$, and $\gamma_u^+$ are random variables. Moreover, since their values are defined for every realization of $w_k$, their expected value with respect to the distribution of $w_k$ can be computed.

We first prove condition \eqref{eq:condru}. Consider the deterministic tightened constraints on the nominal input \eqref{eq:prbbar}, and define the candidate nominal input solution as in (\ref{eq:th1v}), where the dependencies on $x_k$ and $x_{k+1}$ is left implicit to ease the notation, recalling that the realization of the random variable $w_k$ is known at time $k+1$. Moreover, notice that the nominal input law \eqref{eq:th1v} is such that the predicted input $u_{\ell+1|k}^\star$ at time $k$ and the candidate one $u_{\ell|k+1}$  at time $k+1$ are the same, namely
 {
\begin{align}
    u^\star_{\ell+1|k} & = v_{\ell+1|k}^\star + K e_{\ell+1|k}\nonumber\\
    &= v_{\ell+1|k}^\star + K \sum_{i=0}^\ell A_K^{\ell-i} w_{k+i},   \label{eq:u1}\\
u_{\ell|k+1} & = v_{\ell|k+1} + K e_{\ell|k+1}\nonumber\\
&= v_{\ell|k+1} + K \sum_{i=1}^{\ell} A_K^{\ell-i} w_{k+i},\label{eq:u2} 
\end{align}
}
for $\ell \in \N_0^{N-2}$. Indeed, \eqref{eq:u1} and \eqref{eq:u2} are equal if and only if \eqref{eq:th1v} holds. It will be proved next that, since the deterministic constraints \eqref{eq:prbbar} is satisfied by $v_{\ell+1|k}$ with $\gamma_u = \gamma_u^\star$, then it also holds in expectation at $k+1$, i.e., $\E{(v_{\ell|k+1})^\top W_u^{-1}v_{\ell|k+1}} \leq (\gamma_u^\star - \rho (1 - \lambda^{\ell}))^2$.

Let us define the random variable $y = W_u^{-1/2}v_{\ell|k+1}$. From \eqref{eq:th1v}, we get $ y = W_u^{-1/2}(v_{\ell+1|k}^\star + K A_K^{\ell} w_k)$, and consequently $\E{y} = W_u^{-1/2}v_{\ell|k}^\star$ and $\cov{y} = \cov{ W_u^{-1/2} A_K^\ell w_k}$. Knowing that, given a random vector $y$, we have that $\E{y y^\top} = \cov{y} + \E{y}\E{y}^\top$, from \eqref{eq:reachLMI_1} and \eqref{eq:reachLMI_2} we get
\begin{align*}
\E{yy^\top} &= \cov{W_u^{-1/2}K A_K^{\ell} w_k} \\
& \quad + (W_u^{-1/2}v_{\ell|k}^\star)(W_u^{-1/2}v_{\ell|k}^\star)^\top\\
& = \E{W_u^{-1/2} K A_K^{\ell} w_kw_k^\top (A_K^{\ell})^\top K^\top W_u^{-1/2} }\\
& \quad + W_u^{-1/2}v_{\ell+1|k}^\star (v_{\ell+1|k}^\star)^\top W_u^{-1/2}\\
& \preceq \lambda^{2\ell}(1-\lambda)^2  W_u^{-1/2} K W_x K^\top W_u^{-1/2} \\
& \quad + W_u^{-1/2}v_{\ell+1|k}^\star (v_{\ell+1|k}^\star)^\top W_u^{-1/2}\\
& \preceq \lambda^{2\ell}(1-\lambda)^2 \mathbb{I}_n + W_u^{-1/2}v_{\ell+1|k}^\star (v_{\ell+1|k}^\star)^\top W_u^{-1/2}
\end{align*}
where the last inequality follows from the fact that 
\begin{equation*}
 W_u^{-1/2} K W_x K^\top W_u^{-1/2} \preceq \mathbb{I}_n
\end{equation*}
which holds from \eqref{eq:WuWx}. Then, we can compute the expected value of $\E{v_{\ell|k+1}^\top  W_u^{-1}v_{\ell|k+1}}$ relying on the property of random variable for which we have that $\E{y^\top y}=\tr{\E{yy^\top}}$. Hence, we obtain
\begin{align*}
\tr{\E{yy^\top}} & \leq \tr{W_u^{-1/2}v_{\ell+1|k}^\star (v_{\ell+1|k}^\star)^\top W_u^{-1/2}} \\
& \quad +\tr{(1-\lambda)^2 \lambda^{2\ell} \mathbb{I}_n}\\
& = (v_{\ell+1|k}^\star)^\top W_u^{-1} v_{\ell+1|k}^\star + (1-\lambda)^2 \lambda^{2\ell} n\\
& \leq (\gamma_u r_u - \rho (1 - \lambda^{\ell+1}))^2+(1-\lambda)^2 \lambda^{2\ell} n 
\end{align*}
where the last inequality follows from the fact that $v_{\ell+1|k}^\star$ satisfies \eqref{eq:prbbar}. Now, we look for a condition over $\rho$ which guarantees that the constraint \eqref{eq:prbbar} is satisfied in expectation by $v_{\ell|k+1}$ with the same $\gamma_u$, i.e.,
$$\E{v_{\ell|k+1}^\top W_u^{-1} v_{\ell|k+1}}\leq (\gamma_u r_u - \rho (1 - \lambda^{\ell}))^2,$$
which corresponds to define a condition on $\rho$ for which
\begin{equation}\label{eq:th1eta}
(\gamma_u r_u - \rho (1 - \lambda^{\ell+1})^2 + (1-\lambda)^2 \lambda^{2\ell} n\leq (\gamma_u r_u - \rho (1 - \lambda^{\ell}))^2
\end{equation}
holds. Specifically, employing \eqref{eq:constrgamma} and after some manipulation, we can prove that   
\eqref{eq:th1eta} holds if  
\eqref{eq:th1boundrho} is satisfied, which also provides a condition on the violation probability $\phi_\varepsilon(\rho)$.

The analysis above holds for $\ell \in \N_0^{N-2}$, for which the nominal input $v_{\ell|k+1}$ can be defined as a function of $v_{\ell+1|k}^\star$ \eqref{eq:th1v}. Now, let us consider the nominal input defined at time $\ell = N-1$ as in (\ref{eq:th1vN}), given \eqref{eq:th1v} and knowing that $v_{N|k}^\star= K z_{N|k}^\star$.

By defining $y = W_u^{-1/2} v_{N-1|k+1} = W_u^{-1/2} (K z_{N|k} + K A_K^{N-1} w_k)$, the expected value of $\E{y^\top y}$ can be bounded to prove that \eqref{eq:condru} holds also for the last element of the nominal input sequence solving the problem \eqref{eq:SMPC_mixed} at $k+1$. Following a similar procedure as before and applying \eqref{eq:reachLMI_1}, \eqref{eq:reachLMI_2}, \eqref{eq:WuWx}, and \eqref{eq:prcbar2}, we have
\begin{align*}
& \hspace{-0.3cm} \E{yy^\top} = \E{W_u^{-1/2} K A_K^{\ell} w_kw_k^\top (A_K^{\ell})^\top K^\top W_u^{-1/2} }\\
& \quad + W_u^{-1/2} K z_{N|k} z_{N|k}^\top K^\top W_u^{-1/2}\\
& \preceq (1-\lambda)^2 \lambda^{2N-2} \mathbb{I}_n+ W_u^{-1/2}  K z_{N|k} z_{N|k}^\top K^\top W_u^{-1/2}
\end{align*}
from which we get
\begin{equation}\label{eq:th1temp1}
\E{y^\top y} \leq (1-\lambda)^2 \lambda^{2N-2} n + (\gamma_u r_u - \rho (1 - \lambda^{N})^2.
\end{equation}
Hence if the following condition holds
\begin{equation*}
(\gamma_u r_u \! - \rho (1 - \lambda^{N}\!)^2+(1-\lambda)^2 \lambda^{2N-2} n \! \leq \! (\gamma_u r_u - \rho (1 - \lambda^{N-1})^2
\end{equation*}
then the expectation of $v_{N-1|k+1}^\top W_u^{-1}v_{N-1|k+1}$ satisfies the constraint \eqref{eq:prbbar} with the same value of $\gamma_u$, and it can be proved that such condition holds if \eqref{eq:th1boundrho} is satisfied.

Now, let us prove condition \eqref{eq:condrx} related to the constraints on the predicted nominal states. Following an analogous reasoning, it can be proved that $z_{\ell+1|k}$ satisfying \eqref{eq:prabar} implies its satisfaction in expectation also for $z_{\ell|k+1}$. Notice that 
\begin{equation*}
x_{k+1} = z_{0|k+1} = Ax_k + B v_{0|k}^\star + w_k = z_{1|k} + w_k.
\end{equation*}
Moreover, the nominal trajectory at $k+1$, for a given $w_k$ and with nominal control \eqref{eq:th1v} and \eqref{eq:th1vN} is given by 
\begin{align}
z_{1|k+1} & = A z_{0|k+1} + B v_{1|k}^\star + B K w_k \nonumber\\
& = A z_{1|k} + A w_k + B v_{1|k}^\star + B K w_k \nonumber\\
& = z_{2|k} + A_K w_k  \nonumber\\
 \vdots \quad &\nonumber\\
z_{\ell|k+1} & = z_{\ell+1|k} + A_K^{\ell} w_k, \qquad \forall \ell \in \N_0^{N-1}. \label{eq:th1temp2}
\end{align}
Proceeding as above, it can be proved that the random variable $y = W_x^{-1/2}(z_{\ell|k+1}) = W_x^{-1/2}(z_{\ell+1|k} + A_K^\ell w_k)$ is such that 
\begin{equation}\label{eq:th1condz}
\E{y^\top y} \leq (\gamma_x^\star r_x - \rho(1-\lambda^\ell))^2
\end{equation}
if condition \eqref{eq:prabar} holds for $z_{\ell+1|k}$. In fact, since $\E{y} = W_x^{-1/2}z_{\ell+1|k}$ and $\cov{y} = \cov{ W_x^{-1/2} A_K^\ell w_k}$, we have
\begin{align*}
& \hspace{-0.3cm} \E{yy^\top} = \E{W_x^{-1/2} A_K^{\ell} w_kw_k^\top(A_K^{\ell})^\top W_x^{-1/2} }\\
& \quad + W_x^{-1/2}z_{\ell+1|k} (z_{\ell+1|k})^\top W_z^{-1/2}\\
& = W_x^{-1/2} A_K^{\ell} \Gamma_w (A_K^{\ell})^\top W_x^{-1/2} \\
& \quad + W_x^{-1/2}z_{\ell+1|k} (z_{\ell+1|k})^\top W_x^{-1/2}\\
& \preceq (1-\lambda)^2 \lambda^{2\ell} \mathbb{I}_n + W_x^{-1/2}z_{\ell+1|k} (z_{\ell+1|k})^\top W_x^{-1/2},
\end{align*}
from which we get
\begin{align*}
\E{y^\top y} & \leq \tr{(1-\lambda)^2 \lambda^{2\ell} \mathbb{I}_n} \\ 
& \quad + \tr{W_x^{-1/2}z_{\ell+1|k} (z_{\ell+1|k})^\top W_z^{-1/2}} \\
& = (1-\lambda)^2 \lambda^{2\ell} n + (z_{\ell+1|k})^\top W_x^{-1} z_{\ell+1|k} \\
& \leq (1-\lambda)^2 \lambda^{2\ell} n + (\gamma_x r_x - \rho (1 - \lambda^{\ell+1}))^2,
\end{align*}
where the last inequality is holds since $z_{\ell+1|k}$ satisfies \eqref{eq:prabar} with the specific value $\gamma_x = \gamma_x^\star$. Then, using \eqref{eq:constrgamma},  it can be proved that the condition \eqref{eq:th1boundrho} on $\rho$ is sufficient for \eqref{eq:th1condz} to hold.

Finally, let us consider the nominal terminal state at $k+1$, which can be rewritten as follows combining \eqref{eq:th1temp2} and \eqref{eq:th1vN}, i.e.,  
\begin{align*}
z_{N|k+1} & = A z_{N-1|k+1} + B v_{N-1|k+1}\\
& = A (z_{N|k} + A_K^{N-1} w_k) + BK(z_{N|k} + A_K^{N-1}w_k)\\
& = A_K z_{N|k} + A_K^{N}w_k.
\end{align*}
From the same reasoning applied but selecting as random variable $y = W_x^{-1/2}z_{N|k+1}$, one has  
\begin{align}
\E{y^\top y} & \leq (1-\lambda)^2 \lambda^{2N} n + \lambda^2 (\gamma_x r_x - \rho (1 - \lambda^{N}))^2. \label{eq:th1temp3}
\end{align}
Hence, we can prove that \eqref{eq:prcbar} is satisfied in expectation at $k+1$ with the same $\gamma_x$ if 
\begin{equation}\label{eq:th1temp3b}
(1-\lambda)^2 \lambda^{2N} n + \lambda^2 (\gamma_x r_x - \rho (1 - \lambda^{N}))^2 \leq (\gamma_x r_x - \rho (1 - \lambda^{N}))^2
\end{equation}
which holds if 
\begin{equation*}
\lambda^N\left(\sqrt{n\frac{1-\lambda}{1+\lambda}} - \rho \right) + \rho \leq \gamma_x r_x,
\end{equation*}
that is satisfied if \eqref{eq:th1boundrho} holds, which renders the first term non-positive, and since $\rho \leq r_x \leq \gamma_x r_x$.

The last constraint to be proved to hold for the expected solution of problem \eqref{eq:SMPC_mixed} at time $k+1$ is \eqref{eq:prcbar2}. Analogously to  \eqref{eq:th1temp3} and from \eqref{eq:prcbar2} holding with $z_{N|k}$, it follows that \eqref{eq:prcbar2} holds in expectation if 
\begin{equation*}\label{eq:th1temp4}
\lambda^2 (\gamma_u r_u - \rho (1 - \lambda^{N}))^2 + (1-\lambda)^2 \lambda^{2N} n\leq (\gamma_u r_u - \rho (1 - \lambda^{N}))^2,
\end{equation*}
that is satisfied for 
\begin{equation*}
\lambda^N\left(\sqrt{n\frac{1-\lambda}{1+\lambda}} - \rho \right) + \rho \leq \gamma_u r_u.
\end{equation*}
which holds given condition \eqref{eq:th1boundrho} and $\rho \leq r_u \leq \gamma_u r_u$.

\section{Proof of Theorem \ref{th:desc}}
\label{app:desc}
Given and admissible solution of the problem \eqref{eq:SMPC_mixed}, denoted $(\z,\v,\g)$ with $\gamma = (\gamma_x, \gamma_u)$, rewrite the cost of \eqref{eq:SMPC_mixed} in the form $J_N(\mathbf{z},\mathbf{v},\gamma) = J_p(\z,\v) + \eta J_r(\g)$ with $J_p(\cdot)$ and $J_r(\cdot)$ defined in (\ref{eq:Jp}) and (\ref{eq:Jr}), and $J_{LQR}(z) = z^\top P z$ is the optimal cost of the LQR (unconstrained) problem. Moreover, we denote with $J_p^\star(x_k)$ and $J_r^\star(x_k)$ their values at the optimum computed at $x_k$. Finally, we have $(\z^+, \v^+, \g^+)$ as the solution at time $k+1$ obtained by shifting $(\z^\star,\v^\star)$ optimal at time $k$, and $\g^+ = (\g^+_x, \, \g^+_u)$ as defined in Theorem~\ref{th:rplus}. Notice that, from Theorem~\ref{th:rplus}, we have that 
\begin{equation}\label{eq:D2}
\E{J_r(\g^+) | \, w_k } \leq J_r^\star(x_k).
\end{equation}
Moreover, recalling that $x_{k+1}$ and $(\z^+,\v^+,\g^+)$ are deterministic feasible solutions at time $k+1$, given the realization of the disturbance $w_k$ at time $k$, then, given $w_k$, it holds 
\begin{align*}
J_p^\star(x_{k+1}) + \eta J_r^\star(x_{k+1}) \leq & \ J_p(\z^+,\v^+) + \eta J_r(\g^+).
\end{align*}
Taking now the expectation with respect of $w_k$ and from (\ref{eq:D2}), it follows that 
\begin{align}
\E{J_p^\star(x_{k+1}) + \eta J_r^\star(x_{k+1}) \,| \, w_k} \leq \, & \E{J_p(\z^+,\v^+) \,| \, w_k} \nonumber \\
& + \eta J_r^\star(x_k). \label{eq:D3}
\end{align}
Now, to relate the expected value of the optimal cost $J_N^\star(x_{k+1})$ at time $k+1$ conditioned to the realization $w_k$ to $J_N^\star(x_k)$, we use \eqref{eq:D3} and we obtain that 
\begin{align*}
 &\E{J_N^+(\z^+,\v^+,\g^+)\,|\,w_k}-J_N^\star(x_k) = \mathbb {E}\{J_p(\z^+,\v^+)\\ 
 & \hspace{0.5cm} + \eta J_r(\g^+) \,| \, w_k\} - \big(J_p^\star(x_k) + \eta J_r^\star(x_k)\big)\\
 & \hspace{0.5cm} \leq\E{J_p(\z^+,\v^+) \, | \, w_k} - J_p^\star(x_k).
\end{align*}
From the definition of $\z^+$ and $\v^+$ in \eqref{eq:zvplus}, for every $w_k$ one has
\begin{align*}
& \!\!\! J_p(\z^+\!,\v^+) \! = \!\! \sum_{\ell=0}^{N-1} \! \Big(\|z_{\ell|k+1}\|^2_Q \!+\!\|v_{\ell|k+1}\|^2_R\Big)\! +\! \|z_{N|k+1}\|_P^2\\
= & \sum_{\ell=0}^{N-1} \Big(\|z^\star_{\ell+1|k}+A_K^\ell w_k\|^2_Q +\|v^\star_{\ell+1|k}+KA_K^\ell w_k\|^2_R\Big)\\
& + \|A_K z^\star_{N|k} + A_K^{N} w_k\|^2_P\\
= & \sum_{\ell=1}^{N-1} \Big( \|z^\star_{\ell|k} + A_K^{\ell-1} w_k\|^2_Q + \|v^\star_{\ell|k} + K A_K^{\ell-1} w_k\|^2_R \Big)\\
& + \|z^\star_{N|k} + A_K^{N-1} w_k\|^2_Q+ \|K z^\star_{N|k} + K A_K^{N-1} w_k\|^2_R\\
& + \|z^\star_{N|k} + A_K^{N-1} w_k\|^2_{A_K^\top P A_K}\\
= & \sum_{\ell=1}^{N-1} \Big( \|z^\star_{\ell|k} + A_K^{\ell-1} w_k\|^2_Q + \|v^\star_{\ell|k} + K A_K^{\ell-1} w_k\|^2_R \Big)\\
& + \|z^\star_{N|k} + A_K^{N-1} w_k\|^2_P\\
\!\!\! \leq & \sum_{\ell=1}^{N-1} \Big( \|z^\star_{\ell|k}\|^2_Q + \|v^\star_{\ell|k}\|^2_R  \Big) +\|z^\star_{N|k}\|^2_P\\
& + \!\! \sum_{\ell=1}^{N-1}\!\! \Big(\|A_K^{\ell-1} w_k\|^2_Q \!\! + \!\! \|A_K^{\ell-1} w_k\|^2_{K^\top RK}\Big) +\|A_K^{N-1}w_k\|^2_P,
\end{align*}
and consequently
\begin{align*}
& \E{J_p(\z^+,\v^+)\,|\,w_k}
\leq \mathbb{E}\Bigg\{\sum_{\ell=1}^{N-1} \Big( \|z^\star_{\ell|k}\|^2_Q + \|v^\star_{\ell|k}\|^2_R  \Big) \\
& +\|z^\star_{N|k}\|^2_P | w_k \Bigg\} + \mathbb{E}\Bigg\{ \sum_{\ell=1}^{N-1} \Big(\|A_K^{\ell-1} w_k\|^2_Q  \\
& + \|A_K^{\ell-1} w_k\|^2_{K^\top RK}\Big) + \|A_K^{N-1}w_k\|^2_P | w_k \Bigg\}\\
= & \E{ \! \sum_{\ell=1}^{N-1} \! \! \Big( \|z^\star_{\ell|k}\|^2_Q \! + \! \|v^\star_{\ell|k}\|^2_R  \Big) \! + \! \|z^\star_{N|k}\|^2_P + \|w_k\|^2_P \,|\,w_k \! }\\
=& \sum_{\ell=1}^{N-1} \! \! \Big( \|z^\star_{\ell|k}\|^2_Q \! + \! \|v^\star_{\ell|k}\|^2_R  \Big) \! + \! \|z^\star_{N|k}\|^2_P \! + \! \E{\|w_k\|^2_P \,|\,w_k\! }\\
=& \sum_{\ell=1}^{N-1} \Big( \|z^\star_{\ell|k}\|^2_Q + \|v^\star_{\ell|k}\|^2_R  \Big) +\|z^\star_{N|k}\|^2_P + \tr{P\Gamma_w},
\end{align*}
whereas the optimal cost at $k$ is defined as
\begin{align*}
& \hspace{-0.3cm} J^\star_p(x_k) = \sum_{\ell=0}^{N-1} \Big( \|z^\star_{\ell|k}\|^2_Q +\|v^\star_{\ell|k}\|^2_R\big) + \|z^\star_{N|k}\|^2_P\\
= & \|z^\star_{0|k}\|^2_Q \! + \! \|v^\star_{0|k}\|^2_R \! + \! \! \sum_{\ell=1}^{N-1} \! \!  \Big(\|z^\star_{\ell|k}\|^2_Q \! + \! \|v^\star_{\ell|k}\|^2_R \Big) \! + \! \|z^\star_{N|k}\|^2_P\\
= & l(x_k,u^\star_{0|k}) + \sum_{\ell=1}^{N-1} \Big(\|z^\star_{\ell|k}\|^2_Q +\|v^\star_{\ell|k}\|^2_R\big) + \|z^\star_{N|k}\|^2_P,
\end{align*}
since $x_k=z^\star_{0|k}$ and $u^\star_{0|k}=v^\star_{0|k}$.
Hence, we have that
\begin{equation*}
\E{J_p(\z^+,\v^+)\,|\,w_k} - J_p^\star(x_k) \leq \tr{P \Gamma_w} - l(x_k,u^\star_{0|k}).
\end{equation*}
Then, from \eqref{eq:D3}, it follows
\begin{align*}
     & \E{J_p^\star(x_{k+1}) + \eta J_r^\star(x_{k+1}) \,| \, w_k} \leq \E{J_p(\z^+,\v^+) \,| \, w_k} \\ 
     & \! + \eta J_r^\star(x_k) \! \leq \! J_p^\star(x_k) \! + \! \tr{P \Gamma_w} \! -\! l(x_k,u^\star_{0|k}) \! + \eta J_r^\star(x_k),
\end{align*}
from which we obtain an explicit relationship between the optimal cost at time $k$ and the expected value at $k+1$ conditioned to the realization $w_k$, i.e.,
\begin{equation}
\E{J_N^\star(x_{k+1}) \,| \, w_k}-J_N^\star(x_k) \leq \tr{P \Gamma_w}-l(x_k,u^\star_{0|k}).
\label{eq:D7}
\end{equation}
According to \eqref{eq:D7}, the variation of the expected optimal cost over one step, i.e., $\E{J_N^\star(x_{k+1}) \,| \, w_k}-J_N^\star(x_k)$, is strictly lower than zero if $\tr{P\Gamma_w} < l(x_k,u^\star_{0|k})$. Since the expression of $l(x_k,u^\star_{0|k})$ in terms of $x_k$ is in general unknown, though, the region of the state space where the decreasing of $J_N^\star(\cdot)$ is guaranteed in expectation can not be determined explicitly. Nonetheless, since the following   relation holds for every $x_k$ and $u^\star_{0|k}$, i.e.,
\begin{equation}\label{eq:stagecosts}
\|x_k\|^2_Q \leq \|x_k\|^2_Q+\|u^\star_{0|k}\|^2_R = l(x_k,u_{0|k}^\star),
\end{equation}
then a stochastic descent property of the expected value of the optimal solution of the deterministic problem (\ref{eq:SMPC_mixed}) can be established.

Hence, let us first consider the geometric drift condition result presented in \cite[Proposition 1]{chatterjee2014stability}, claiming that if $\E{J_N^\star(x_{k+1})} \leq (1-\mu) J_N^\star(x_k)$ for all $x_k \notin \Omega$, with $\Omega$ compact and $\mu \in (0, 1)$, then $\E{J_N^\star(x_{k+j})} \leq (1-\mu)^jJ_N^\star(x_k) + \beta/\mu$ with $\beta = \sup_{x \in  \Omega} \E{J_N(x)}$. \MF{Moreover, let us consider the function defined as $\hat{J}_N(x) = J_{LQR}(x) + \eta J_r(\hat{\gamma})$, with $\hat{\gamma} = (\hat{\gamma}_x, \hat{\gamma}_u)$ given by 
\begin{align*}
\hat{\gamma}_x \!\! = \! \max \! \left(\!\!\frac{\sqrt{x^\top W_x^{-1} x}}{r_x}, 1 \!\right)\!\!, \,
\hat{\gamma}_u \!\! = \! \max \! \left(\!\!\frac{\sqrt{x^\top W_x^{-1} x}}{r_u}, 1 \!\right)\!,
\end{align*}
which is, in practice, the cost of problem \eqref{eq:SMPC_mixed} if the LQR solution is applied and $\hat{\gamma}_x r_x = \hat{\gamma}_u r_u \geq \sqrt{x^\top W_x^{-1} x}$, from Proposition~\ref{prop:inEcWx} with $\sqrt{x^\top W_x^{-1} x}$ in spite of $r_{xu}$.}  Hence, from Proposition~\ref{prop:inEcWx} we have
\begin{equation}\label{eq:relationVs}
J_{LQR}(x_k) \leq J_N^\star(x_k) \leq \hat{J}_N(x_k), \quad \forall x_k \in \R^n 
\end{equation} 
and 
\begin{equation}\label{eq:relationVsegal}
J_{LQR}(x_k) = J_N^\star(x_k) = \hat{J}_N(x_k), \quad \forall x_k \in \Ec_{W_x} (r_{xu}).
\end{equation} 
Then, from (\ref{eq:stagecosts})--(\ref{eq:relationVs}) if $x_k \notin \mathbb{X}_\mu$ with 
\begin{equation*}
\mathbb{X}_\mu = \{x_k \in \X \,| \ \tr{P \Gamma_w} - x_k^\top Q x_k \geq - \mu \hat{J}_N(x_k) \},
\end{equation*}
we obtain that 
\begin{equation}
\tr{P \Gamma_w} - l(x_k,u_{0|k}^\star) \leq - \mu J_N^\star(x_k),
\end{equation}
and consequently the geometric drift condition holds for all $x_k \notin \X_\mu$. Moreover, if the inclusion $\mathbb{X}_\mu \subseteq \Ec_{W_x} (r_{xu})$ holds, implied by (\ref{eq:LMIXa}), then, from (\ref{eq:relationVsegal}), the set $\mathbb{X}_\mu$ is given by the following ellipsoid
\begin{equation*}
\mathbb{X}_\mu = \{x_k \in \X \,| \ x_k^\top (Q - \mu P) x_k \leq \tr{P \Gamma_w}\}.
\end{equation*}
Therefore, if (\ref{eq:LMIXa}) holds, then condition (\ref{eq:convergence}) is satisfied for all $i \in \N$ with $\beta$ defined as (\ref{eq:beta}). Hence, $\E{J_N^\star(x_{k+i})}$ converges to $\beta/\mu$ in expectation. 

\MF{Finally, note that condition (\ref{eq:LMIXa}), together with (\ref{eq:relationVs}) and (\ref{eq:relationVsegal}), implies also that the value $\beta/\mu$, toward which the MPC cost converges in expectation from (\ref{eq:convergence}), is in the interior of $\Ec_{W_x}(r_{xu})$. Therefore, the state converges within the set $\Ec_{W_x}(r_{xu})$ and the relaxations eventually vanish, i.e., (\ref{eq:limitr}) hold. }

\bibliographystyle{IEEEtran}
\bibliography{SMPC}

\begin{thebibliography}{10}
\providecommand{\url}[1]{#1}
\csname url@rmstyle\endcsname
\providecommand{\newblock}{\relax}
\providecommand{\bibinfo}[2]{#2}
\providecommand\BIBentrySTDinterwordspacing{\spaceskip=0pt\relax}
\providecommand\BIBentryALTinterwordstretchfactor{4}
\providecommand\BIBentryALTinterwordspacing{\spaceskip=\fontdimen2\font plus
\BIBentryALTinterwordstretchfactor\fontdimen3\font minus
  \fontdimen4\font\relax}
\providecommand\BIBforeignlanguage[2]{{%
\expandafter\ifx\csname l@#1\endcsname\relax
\typeout{** WARNING: IEEEtran.bst: No hyphenation pattern has been}%
\typeout{** loaded for the language `#1'. Using the pattern for}%
\typeout{** the default language instead.}%
\else
\language=\csname l@#1\endcsname
\fi
#2}}

\bibitem{Prekopa1995}
A.~Pr\'ekopa, \emph{Stochastic Programming}.\hskip 1em plus 0.5em minus
  0.4em\relax Kluwer Academic Publishers, 1995.

\bibitem{mayne2005robust}
D.~Q. Mayne, M.~M. Seron, and S.~V. Rakovi{\'c}, ``Robust model predictive
  control of constrained linear systems with bounded disturbances,''
  \emph{Automatica}, vol.~41, no.~2, pp. 219--224, 2005.

\bibitem{farina2016stochastic}
M.~Farina, L.~Giulioni, and R.~Scattolini, ``Stochastic linear model predictive
  control with chance constraints -- {A} review,'' \emph{Journal of Process
  Control}, vol.~44, pp. 53--67, 2016.

\bibitem{lorenzen2016constraint}
M.~Lorenzen, F.~Dabbene, R.~Tempo, and F.~Allg{\"o}wer, ``Constraint-tightening
  and stability in stochastic model predictive control,'' \emph{IEEE
  Transactions on Automatic Control}, vol.~62, no.~7, pp. 3165--3177, 2016.

\bibitem{munoz2020convergence}
D.~Munoz-Carpintero and M.~Cannon, ``Convergence of stochastic nonlinear
  systems and implications for stochastic model-predictive control,''
  \emph{IEEE Transactions on Automatic Control}, vol.~66, no.~6, pp.
  2832--2839, 2020.

\bibitem{Blackmore}
L.~Blackmore, M.~Ono, A.~Bektassov, and B.~C. Williams, ``A probabilistic
  particle-control approximation of chance-constrained stochastic predictive
  control,'' \emph{IEEE Transactions on Robotics}, vol.~26, no.~3, pp.
  502--517, 2010.

\bibitem{CalafioreFagiano}
G.~C. Calafiore and L.~Fagiano, ``Robust model predictive control via scenario
  optimization,'' \emph{IEEE Transactions on Automatic Control}, vol.~58,
  no.~1, pp. 219--224, 2013.

\bibitem{hewing2018stochastic}
L.~Hewing and M.~N. Zeilinger, ``Stochastic model predictive control for linear
  systems using probabilistic reachable sets,'' in \emph{2018 IEEE Conference
  on Decision and Control (CDC)}.\hskip 1em plus 0.5em minus 0.4em\relax IEEE,
  2018, pp. 5182--5188.

\bibitem{kouvaritakis2016model}
B.~Kouvaritakis and M.~Cannon, \emph{Model Predictive Control: Classical,
  Robust, and Stochastic}.\hskip 1em plus 0.5em minus 0.4em\relax Springer,
  2015.

\bibitem{rawlings2017model}
J.~B. Rawlings, D.~Q. Mayne, and M.~Diehl, \emph{Model predictive control:
  {T}heory, computation, and design}.\hskip 1em plus 0.5em minus 0.4em\relax
  Nob Hill Publishing Madison, WI, 2017, vol.~2.

\bibitem{kouvaritakis2010explicit}
B.~Kouvaritakis, M.~Cannon, S.~V. Rakovic, and Q.~Cheng, ``Explicit use of
  probabilistic distributions in linear predictive control,''
  \emph{Automatica}, vol.~46, no.~10, pp. 1719--1724, 2010.

\bibitem{goulart2006optimization}
P.~J. Goulart, E.~C. Kerrigan, and J.~M. Maciejowski, ``Optimization over state
  feedback policies for robust control with constraints,'' \emph{Automatica},
  vol.~42, no.~4, pp. 523--533, 2006.

\bibitem{paulson2020stochastic}
J.~A. Paulson, E.~A. Buehler, R.~D. Braatz, and A.~Mesbah, ``Stochastic model
  predictive control with joint chance constraints,'' \emph{International
  Journal of Control}, vol.~93, no.~1, pp. 126--139, 2020.

\bibitem{farina2013probabilistic}
M.~Farina, L.~Giulioni, L.~Magni, and R.~Scattolini, ``A probabilistic approach
  to model predictive control,'' in \emph{52nd IEEE Conference on Decision and
  Control {(CDC)}}.\hskip 1em plus 0.5em minus 0.4em\relax IEEE, 2013, pp.
  7734--7739.

\bibitem{di1994exact}
G.~Di~Pillo, ``Exact penalty methods,'' \emph{Algorithms for continuous
  optimization: the state of the art}, pp. 209--253, 1994.

\bibitem{yan2018stochastic}
S.~Yan, P.~Goulart, and M.~Cannon, ``Stochastic model predictive control with
  discounted probabilistic constraints,'' in \emph{2018 European Control
  Conference (ECC)}.\hskip 1em plus 0.5em minus 0.4em\relax IEEE, 2018, pp.
  1003--1008.

\bibitem{mammarella2020probabilistic}
M.~Mammarella, T.~Alamo, S.~Lucia, and F.~Dabbene, ``A probabilistic validation
  approach for penalty function design in stochastic model predictive
  control,'' \emph{IFAC-PapersOnLine}, vol.~53, no.~2, pp. 11\,271--11\,276,
  2020.

\bibitem{TeBaDa:97}
R.~Tempo, E.~Bai, and F.~Dabbene, ``Probabilistic robustness analysis: explicit
  bounds for the minimum number of samples,'' \emph{Systems {\&} Control
  Letters}, vol.~30, pp. 237--242, 1997.

\bibitem{Alamo:15}
T.~Alamo, R.~Tempo, A.~Luque, and D.~Ramirez, ``Randomized methods for design
  of uncertain systems: sample complexity and sequential algorithms,''
  \emph{Automatica}, vol.~52, pp. 160--172, 2015.

\bibitem{kerrigan2000soft}
E.~C. Kerrigan and J.~M. Maciejowski, ``Soft constraints and exact penalty
  functions in model predictive control,'' in \emph{Control 2000 Conference},
  2000, pp. 2319--2327.

\bibitem{karg2019probabilistic}
B.~Karg, T.~Alamo, and S.~Lucia, ``Probabilistic performance validation of deep
  learning-based robust nmpc controllers,'' \emph{International Journal of
  Robust and Nonlinear Control}, vol.~31, no.~18, pp. 8855--8876, 2021.

\bibitem{farina2015approach}
M.~Farina, L.~Giulioni, L.~Magni, and R.~Scattolini, ``An approach to
  output-feedback mpc of stochastic linear discrete-time systems,''
  \emph{Automatica}, vol.~55, pp. 140--149, 2015.

\bibitem{hewing2020recursively}
L.~Hewing, K.~P. Wabersich, and M.~N. Zeilinger, ``Recursively feasible
  stochastic model predictive control using indirect feedback,''
  \emph{Automatica}, vol. 119, p. 109095, 2020.

\bibitem{fiacchini2021probabilistic}
M.~Fiacchini and T.~Alamo, ``Probabilistic reachable and invariant sets for
  linear systems with correlated disturbance,'' \emph{Automatica}, vol. 132, p.
  109808, 2021.

\bibitem{kohler2022recursively}
J.~K{\"o}hler and M.~N. Zeilinger, ``Recursively feasible stochastic predictive
  control using an interpolating initial state constraint,'' \emph{IEEE Control
  Systems Letters}, vol.~6, pp. 2743--2748, 2022.

\bibitem{schluter2022stochastic}
H.~Schl{\"u}ter and F.~Allg{\"o}wer, ``Stochastic model predictive control
  using initial state optimization,'' \emph{IFAC-PapersOnLine}, vol.~55,
  no.~30, pp. 454--459, 2022.

\bibitem{schluter2023stochastic}
------, ``Stochastic model predictive control using initial state and variance
  interpolation,'' in \emph{2023 IEEE Conference on Decision and Control
  (CDC)}.\hskip 1em plus 0.5em minus 0.4em\relax IEEE, 2023, pp. 6700--6706.

\bibitem{mammarella2018offline}
M.~Mammarella, M.~Lorenzen, E.~Capello, H.~Park, F.~Dabbene, G.~Guglieri,
  M.~Romano, and F.~Allg{\"o}wer, ``An offline-sampling {SMPC} framework with
  application to autonomous space maneuvers,'' \emph{IEEE Transactions on
  Control Systems Technology}, vol.~28, no.~2, pp. 388--402, 2018.

\bibitem{kofman2012probabilistic}
E.~Kofman, J.~A. De~Don{\'a}, and M.~M. Seron, ``Probabilistic set invariance
  and ultimate boundedness,'' \emph{Automatica}, vol.~48, no.~10, pp.
  2670--2676, 2012.

\bibitem{hewing2018correspondence}
L.~Hewing, A.~Carron, K.~P. Wabersich, and M.~N. Zeilinger, ``On a
  correspondence between probabilistic and robust invariant sets for linear
  systems,'' in \emph{2018 European Control Conference (ECC)}.\hskip 1em plus
  0.5em minus 0.4em\relax IEEE, 2018, pp. 1642--1647.

\bibitem{arcari2023stochastic}
E.~Arcari, A.~Iannelli, A.~Carron, and M.~N. Zeilinger, ``Stochastic {MPC} with
  robustness to bounded parametric uncertainty,'' \emph{IEEE Transactions on
  Automatic Control}, vol.~68, no.~12, pp. 7601--7615, 2023.

\bibitem{Lew2019ChanceConstrainedOA}
T.~Lew, F.~Lyck, and G.~M{\"u}ller, ``Chance-constrained optimal altitude
  control of a rocket,'' in \emph{8th European Conference for Aeronautics and
  Space Sciences (EUCASS)}, 2019, pp. 1--15.

\bibitem{Gross2021mission}
K.~Wang and S.~Gros, ``Recursive feasibility of stochastic model predictive
  control with mission-wide probabilistic constraints,'' in \emph{2021 IEEE
  Conference on Decision and Control (CDC)}, 2021, pp. 2312--2317.

\bibitem{raghuraman2022long}
V.~Raghuraman and J.~P. Koeln, ``{Long duration stochastic MPC with
  mission-wide probabilistic constraints using waysets},'' \emph{IEEE Control
  Systems Letters}, vol.~7, pp. 865--870, 2022.

\bibitem{chatterjee2014stability}
D.~Chatterjee and J.~Lygeros, ``On stability and performance of stochastic
  predictive control techniques,'' \emph{IEEE Transactions on Automatic
  Control}, vol.~60, no.~2, pp. 509--514, 2014.

\end{thebibliography}

\end{document}